\documentclass[11pt,oneside,reqno]{amsart}
\usepackage[utf8]{inputenc}
\usepackage[T1]{fontenc}
\usepackage{lmodern}
\usepackage{amssymb, amsmath, amsthm, mathtools, cases,enumerate}
\usepackage{color}
\usepackage{hyperref,mathscinet}
\usepackage{natbib}
\usepackage{microtype}
\usepackage{hyphenat}

\usepackage{theoremref}

\newtheorem{theorem}{Theorem}[section]
\newtheorem{lemma}[theorem]{Lemma}
\newtheorem{proposition}[theorem]{Proposition}

\theoremstyle{definition}

\newtheorem{remark}[theorem]{Remark}
\theoremstyle{remark}

\numberwithin{equation}{section}

\allowdisplaybreaks

\newcommand{\C}{\mathcal{C}}

\begin{document}
\setcitestyle{numbers}
\bibliographystyle{plainnat}

\setcounter{page}{1}

\begin{abstract}

It is well known that the optimal Sobolev embedding (i.e. embedding of Sobolev space into the smallest possible target space) is non-compact. Recently it was proved that such embedding is maximally non-compact and even not strictly singular (i.e. there exists an infinitely dimensional subspace on which the embedding is invertible). In this paper, we investigate whether strict singularity is a global or localized property. We construct Sobolev embedding which is locally compact at the neighborhood of every point within their domain, except one single point. For this embedding, we obtain a sharp condition that distinguishes compactness from non-compactness and observe that in the context of Sobolev embeddings, non-compactness occurring just at only one point within the domain could give rise to non-compact mapping that is not strictly singular. Furthermore, we establish lower bounds for the Bernstein numbers, entropy numbers, and the measure of non-compactness.

\end{abstract}

\title[]
{Quality of non-compactness for Sobolev Embedding with one point non-compactness.}

\author{Chian Yeong Chuah, Jan Lang}

\address{Chian Yeong Chuah, Department of Mathematics, The Ohio State University, 231 W. 18th Ave.
Columbus, OH 43210, USA.}
\email{\textcolor[rgb]{0.00,0.00,0.84}{chuah.21@osu.edu}}

\address{Jan Lang, Department of Mathematics, The Ohio State University, 231 W. 18th Ave.
Columbus, OH 43210, USA, {\it and }
Czech Technical University in Prague, Faculty of Electrical Engineering, Department of Mathematics, Technická 2, 166 27 Praha 6, Czech Republic}

\email{\textcolor[rgb]{0.00,0.00,0.84}{lang.162@osu.edu}}


\maketitle

\section{Introduction}

Embeddings of Sobolev spaces into Lebesgue or Lorentz spaces, defined on an open set $\Omega \subset \mathbb{R}^d$, hold significant importance in both the theory of partial differential equations and approximation theory. Understanding the ``quality'' and internal characteristics of these embeddings is essential for numerous practical applications.

Let us consider the following two Sobolev embeddings:
\begin{equation} \label{1}
	I: W_0^{1,p}(\Omega) \to L^{p^*}(\Omega)
\end{equation}
and 
\begin{equation} \label{2}
	I: W_0^{1,p}(\Omega) \to L^{p^*,p}(\Omega),
\end{equation}
where $p\in [1,d)$ and $p^*=dp/(d-p)$. By $\Omega$ we denote a bounded open subset of $\mathbb{R}^d$ which is sufficiently regular (for example, Lipschitz), and $W_0^{1,p}(\Omega)$ refers to the Banach space of  all functions from the Sobolev space $W^{1,p}(\Omega)$ with zero boundary values (see the next section for definitions).

Both embeddings (\ref{1}) and (\ref{2}) are non-compact and the target spaces are optimal. The Lebesgue space $L^{p^*}(\Omega)$ is the optimal target space in (\ref{1}) among all Lebesgue spaces - i.e. for any smaller Lebesgue space $L^q(\Omega) \subsetneqq L^{p^*}(\Omega)$ ($p^* <q $), the embedding $I:W_0^{1,p}(\Omega) \to L^{q}(\Omega)$ is unbounded. It is known  (see \cite{Pe}) that (\ref{1}) can be improved if one looks for the optimal target space among all Lorentz spaces. In this case, the Lorentz space $L^{p^*,p}(\Omega)$ is the optimal target space. Since $L^{p^*,p}\subsetneqq L^{p^*}$, (\ref{2}) is an improvement over (\ref{1}).
In addition, the  Lorentz space $L^{p^*,p}(\Omega)$ is the optimal target space  among all the rearrangement-invariant function spaces (i.e. collections of function spaces which include, among others,  Lebesgue spaces, Lorentz spaces, Orlicz spaces, and Marcinkiewicz spaces). This means that if $I:W_0^{1,p}(\Omega) \to Y(\Omega)$ is valid, then $L^{p^*,p}(\Omega) \subseteqq Y(\Omega)$ (see \cite{KP2006}).

It is worth mentioning that both embeddings (\ref{1}) and (\ref{2}) are not only non-compact, but they are also ``maximally non-compact'' as their norms are equal to their measures of non-compactness (see [Definition 2.7]\cite{EdEv}). This was proved in \cite{He} and \cite{Bo2020}. In \cite{LMOP}, it was showed that when the target space $L^{p^*,p}(\Omega)$ in (\ref{2}) is enlarged to the larger Lebesgue space  $L^{p^*,r}(\Omega)$ $(p< r\le \infty)$, for which we have $L^{p^*,p} \subsetneqq L^{p^*,r}(\Omega) \subset L^{p^*, \infty}(\Omega)$, then  the corresponding embedding is still non-compact but surprisingly  also ``maximally non-compact''.

Given strict difference between $L^{p^*,r}$ spaces with different $r$, one would expect that the quality of the corresponding Sobolev embedding into $L^{p^*,r}$ spaces should depend on $r$.  However, as previously noted, these differences are not fully captured by the measure of non-compactness alone. This prompts the need for a more refined approach.

The concepts of strict singularity and Bernstein numbers emerge as promising candidates for quantifying the ``quality'' of non-compactness. By delving into these quantities, we can gain deeper insights into the nuanced inner structures of these Sobolev embeddings.
 This was confirmed by  \cite{LaMi}, in which it was shown that (\ref{2}) is not strictly singular (i.e. there exists an infinite dimensional subspace on which the embedding is isomorphic) and that the Bernstein numbers of (\ref{2}) are equal to the norm of the embedding, and that (\ref{1}) is finitely strictly singular (i.e. Bernstein numbers converge to 0). The main techniques used in  \cite{LaMi} were based on the existence of non-compactness at each point, i.e. the embedding has the same measure of non-compactness at each open ball in the underlying domain $\Omega$, and  the norm of embedding is invariant with respect to a re-scaling.
 
 The natural question that arises from  \cite{LaMi},  \cite{LMOP}, \cite{He}, and \cite{Bo2020} is: which ``quality'' of non-compactness can be expected in the case when the norm of embedding is not invariant with respect to re-scaling and the embedding is compact at each neighborhood except at a neighborhood of just one point in the domain. 
 
 In this paper, we construct a Sobolev embedding from $W^{1,p}_0(\Omega)$ into variable Lorentz space $L^{q(\cdot),p}(\Omega)$ which is  non-compact in just a neighborhood of one point but compact in all other neighborhoods in $\Omega$. Moreover, we show that this Sobolev embedding is not strictly singular (\thref{mainTh}), and the Bernstein numbers and the measure of non-compactness have strictly positive lower bound (\thref{mainresult}). By this, we demonstrate that a localized non-compactness at just one point is satisfactory to generate an infinitely dimensional system of functions on which the corresponding Sobolev embedding is invertible (i.e. Sobolev embedding is not strictly singular). Consequently, we obtain an embedding which is essentially as non-compact as the most non-compact case, i.e. embedding (\ref{2}).   
 
 As a by-product of our techniques, we extend results from \cite{EGN} onto variable Lorentz spaces and we obtain sharp conditions which guarantee compactness and non-compactness and we also extend results from \cite{KS2008} into variable Lorentz spaces.      
 
The paper is structured as follows. In the next section, we recall some basic definitions and notations. In Section 3, we start by introducing almost compact embeddings, make observations about almost compact embeddings in the context of variable Lorentz spaces (\thref{almostcompact}) and describe conditions under which our embeddings are compact (\thref{q-compact}) or non-compact (\thref{propnon-compact}). In the last section, we obtain the main results (\thref{mainresult}) and (\thref{mainTh}).   

\section{Preliminaries}

We start this section by recalling definitions of function spaces that  will be used throughout this paper.

Given a measurable function $f : \Omega \to \mathbb{C}$, the distribution function of $f$, $d_{f} : [0, \infty) \to [0, \infty)$ is defined by:
\begin{displaymath}
    d_{f} (\lambda) := \left| \left\{ x \in \Omega : |f(x)| > \lambda \right\} \right|.
\end{displaymath}
The decreasing rearrangement of $f$, $f^{*} : [0, \infty) \to [0, \infty)$ is defined by:
\begin{displaymath}
    f^{*} (t) := \inf \left\{ s > 0 : d_{f} (s) \leq t \right\}.
\end{displaymath}

Let $A \subseteq \mathbb{R}^{d}$ be a Lebesgue measurable set. The symmetric decreasing rearrangement of $A$, $A^{\#}$ is defined by:
\begin{displaymath}
    A^{\#} := \left\{ x \in \mathbb{R}^{d} : \omega_{d} |x|^{d} < \left| A \right| \right\}, \text{ where } \omega_{d} \text{ is the volume of the } d \text{ dimensional unit ball.}
\end{displaymath}
Clearly, $\left| A^{\#} \right| = \left| A \right|$.

Given a measurable function $f : \Omega \to \mathbb{C}$, the symmetric decreasing rearrangement of $f$, $f^{\#} : \mathbb{R}^{d} \to [0, \infty)$ is defined by:
\begin{displaymath}
    f^{\#} (x) := \int_{0}^{\infty} \chi_{\left\{ y \in \Omega : | f(y) | > t \right\}^{\#} } (x) d t.
\end{displaymath}

Let $p \in [1, \infty)$ and $q : \Omega \to [1, \infty]$ be a measurable function. The number $q_{+} (\Omega)$ is defined as $ q_{+} (\Omega) := \sup_{x \in \Omega} q(x),$ and $\Omega_{\infty} := \left\{ x \in \Omega : q (x) = \infty \right\}$. The modular functional associated with $q(\cdot)$, $\rho_{q(\cdot), \Omega}$ is defined as:

\begin{displaymath}
    \rho_{q(\cdot), \Omega} (f) := \int_{\Omega \setminus \Omega_{\infty}} | f(x) |^{q(x)} d x + \left\| f \right\|_{L^{\infty} (\Omega_{\infty})}.
\end{displaymath} 
By $L^{q(\cdot)}(\Omega)$ and $L^{q(\cdot),p }(\Omega)$, we denote the variable exponent Lebesgue space and the variable exponent Lorentz space, respectively, which are defined via the following norms:

\begin{displaymath}
    \left\| f \right\|_{L^{q(\cdot)} (\Omega)} := \inf \left\{ \lambda > 0 : \rho_{q(\cdot), \Omega} \left( \frac{f}{\lambda} \right) \leq 1 \right\} = \sup \left\{ \lambda > 0 : \rho_{q(\cdot), \Omega} \left( \frac{f}{\lambda} \right) > 1 \right\},
\end{displaymath}

\begin{displaymath}
\left\| f \right\|_{L^{q(\cdot), p} (\Omega)} := \left[ \int_{0}^{\infty} \lambda^{p - 1} \left\| \chi_{\left\{ x \in \Omega : |f(x)| > \lambda \right\}} \right\|_{L^{q(\cdot)} (\Omega)}^{p} \ d \lambda \right]^{\frac{1}{p}}.   
\end{displaymath}
These spaces can be considered as generalizations of the "standard" Lebesgue or Lorentz spaces when $q(\cdot)$ is replaced by a constant. More information about these spaces can be found in \cite{KovRak} and \cite{KeVy}.

The set of all smooth (i.e., infinitely differentiable) functions that are compactly supported in $\Omega$ is denoted by $\C_0^{\infty}(\Omega)$.

By $W^{1,p}(\Omega)$ we denote {\it the classical first-order Sobolev space} on $\Omega$ defined by the norm
\begin{displaymath}
\left\| f \right\|_{W^{1, p} (\Omega)} := \left[ \int_{\Omega}|f(x)|^p d x + \int_{\Omega} | (\nabla f) (x) |^{p} \ d x \right]^{\frac{1}{p}}. 
\end{displaymath}

We denote the closure of $\C_0^{\infty}(\Omega)$ in $W^{1,p}(\Omega)$ by $W_0^{1,p}(\Omega)$ and equip it with the norm
\begin{displaymath}
	\left\| f \right\|_{W_0^{1, p} (\Omega)} := \left[  \int_{\Omega} | (\nabla f) (x) |^{p} \ d x \right]^{\frac{1}{p}},
\end{displaymath}
which is equivalent with the $W^{1,p}$ norm on $\C_0^{\infty}(\Omega).$



Now, we recall the definitions of some s-numbers (Bernstein and Kolmogorov numbers) and entropy  numbers.

Let $X$, $Y$ be Banach spaces and $T \in B(X, Y)$. We denote  $BX=\{ x \in X : \|x\|_X \le 1\}$.

The $n$-th  Kolmogorov number is defined by
\begin{equation*}
	d_n(T):= \inf_{Y_n}  \sup_{y\in T(BX)} \inf_{z\in Y_n} \|y-z\|_Y,
\end{equation*} 
where the infimum is taken over all $n$-dimensional sub-spaces $Y_n$ of $Y$ and $T(BX)=\{T(x) : x\in BX\}$.

The $n$-th Bernstein number of $T$, $b_{n} (T)$ is defined by:
\begin{displaymath}
	b_{n} (T) := \sup \left\{ \inf_{x \in X_{n}, \| x \|_{X} = 1} \left\| T(x) \right\|_{Y} : X_{n} \text{ is an } n \text{ dimensional subspace of } X \right\}.
\end{displaymath}

The $s$-numbers defined above (see \cite{EdEv}) are ordered $\|T\| \ge d_n(T) \ge b_n(T)$ and  can be used for describing the ``quality''  of operators.  It is worth noting that $\alpha(T):=\lim_{n\to \infty} d_n(T)=0$ if and only if $T$ is compact.

Another quantity that is used for describing compactness are entropy numbers. We say that the $n$-the entropy number of $T$ is defined by
\[ e_n(T):= \inf \{\varepsilon>0: T(BX) \mbox{ can be covered by }2^{n-1} \mbox{ balls in }Y \mbox{ with radius }\varepsilon  \}.
\]

Note that $T$ is compact if and only if $\beta(T):= \lim_{n \to \infty} e_n(T) =0$, where $\beta(T)$ is called {\it the measure of non-compactness}. If $\beta(T)=\|T\|$, we say that $T$ is {\it maximally non-compact}. 

We say that an operator $T\in B(X,Y)$ is {\it strictly singular} (SS) if there is no infinite-dimensional closed subspace $Z$ of $X$ such that the restriction $T|_Z$ is an isomorphism of $Z$ onto $T(Z)$. Equivalently, this can be  described that for each infinite dimensional (closed) subspace $Z$ of $X$, we have
\[ \inf \{ \|T (x) \|_Y: \|x\|_X=1, x\in Z\}=0.
\]  

An operator $T\in B(X,Y)$ is said to be {\it finitely strictly singular} (FSS) if for any
 $\varepsilon >0$, there exists $N(\varepsilon) \in \mathbb{N}$ such that if $E$ is a subspace of $X$ with $\dim E \ge N(\varepsilon)$, 
 then there exists $x\in E, \|x\|_X=1,$ such that $\|T(x) \|_Y \le \varepsilon$.
 
It is possible to see that an operator $T$ is finitely strictly singular if and only if
\[ \lim_{n \to \infty} b_n(T)=0.
\]
We conclude this section with the well-known fact:
\[Compact \subset FSS \subset SS.\]
  
\section{Compactness and Non-compactness conditions}

In this section, we study conditions on $q(\cdot)$
under which the Sobolev embedding  \[I:W_0^{1, p} (\Omega) \to L^{q(\cdot), p} (\Omega) \] is compact or non-compact. We consider the case when $1\le q(\cdot) \le p^*$ and $q(\cdot)$ is approaching $p^*$ only at one point of the domain $x_0 \in \Omega$. This will bring us to the situation in which $I$ is ``locally'' compact at the neighborhood of each point in $\Omega \setminus \{x_0\}$ and possibly non-compact  only at neighborhoods of $x_0$.

In order to obtain conditions on the growth of $q(\cdot)$ at $x_0$ which will guarantee compactness/non-compactness, we need to make some observations.  We start with statements needed for describing almost compact embeddings on variable Lorentz spaces  which is  a generalization of the main results from \cite{EGN}. Employing these results for variable Lorentz spaces gives us  a condition  on $q(\cdot)$ which guarantees the compactness of Sobolev embedding. 

First, we recall some definitions and the known results: If $X(\Omega)$ and $Y(\Omega)$ are Banach function spaces over $\Omega$, then we say that $X(\Omega)$ is {\it almost compactly embedded} into $Y(\Omega)$, denoted by $X(\Omega) \overset{*}\hookrightarrow Y(\Omega)$, if for every sequence $(E_n)_{n = 1}^{\infty}$ of measurable subsets of $\Omega$ such that $E_n \to 0$ a.e., we have
\[\lim_{n\to \infty} \sup_{\|u\|_X\le 1} \|u \chi_{E_n} \|_Y = 0. \]

The next proposition demonstrates the connection with the compactness of Sobolev embeddings.
\begin{proposition} \label{prop3.4}
    Let $X, Y,$ and $Z$ be Banach function spaces over bounded set $\Omega$ and assume
    \[
    W^1(X) \hookrightarrow Y, \mbox{ and } \quad Y \overset{*}\hookrightarrow Z,
    \]
    where $W^1(X)$ denotes the Sobolev space with norm $\|u\|_{W^{1} (X)} = \|u\|_X +\|\nabla u\|_X.$
    
    Then $W^1(X) \hookrightarrow Z$ is a compact embedding.
\end{proposition}
\begin{proof}
    See \cite[Prop 2.2]{EGN} or literature noted there.
\end{proof}

\begin{lemma} \thlabel{exponent}
    Let $s : \Omega \to \mathbb{R}$ be a measurable function and $\alpha > 1$. Then, for each $t > 0$, $\left[ \alpha^{s(\cdot)} \right]^{*} (t) = \alpha^{s^{*}(t)}$. 
\end{lemma}
\begin{proof}
This is proved in \cite[Lemma 2.10]{EGN}.
\end{proof}

\begin{lemma} \label{3.5}
     Let $p : \Omega \to [1, \infty)$ and $q : \Omega \to [1, \infty)$ be measurable functions where $q(x) \leq p(x) \leq p_{+} < \infty$ for all $x \in \Omega$. Assume that $\left| \left\{ x \in \Omega : p(x) = q(x) \right\} \right| = 0$ and for any sequence $(E_{n})_{n = 1}^{\infty}$ of measurable subsets of $\Omega$ such that $| E_{n}| \to 0$, we have 

    \begin{displaymath}
        \left\| \chi_{E_{n}} \right\|_{L^{\frac{p(\cdot) q(\cdot)}{p(\cdot) - q(\cdot)}} (\Omega)} \to 0,
    \end{displaymath}
where $p(\cdot)q(\cdot)/(p(\cdot)-q(\cdot))=\infty$ when $p(\cdot)=q(\cdot).$
    Then, $L^{p(\cdot), t} (\Omega)$ is almost compactly embedded into $L^{q(\cdot), t} (\Omega)$ for any $t \in [1, \infty)$ .
\end{lemma}

\begin{proof}
    Let $u \in L^{p(\cdot), t} (\Omega)$ where $\left\| u \right\|_{L^{p(\cdot), t} (\Omega)} \leq 1$. Let $(E_{n})_{n = 1}^{\infty}$ be a sequence of measurable subsets of $\Omega$ where $\chi_{E_{n}} \to 0$ a.e. as $n \to \infty$.  
    
    Now, fix $n \in \mathbb{N}$ and $\lambda > 0$. By Holder's inequality, there exists $K \in (0, \infty)$ such that 
    \begin{align*}
        & \left\| \chi_{\left\{ x \in \Omega : \left| u (x) \right| > \lambda \right\}} \cdot \chi_{E_{n}} \right\|_{L^{q(\cdot)} (\Omega)} \leq K \left\| \chi_{\left\{ x \in \Omega : \left| u (x) \right| > \lambda \right\}} \right\|_{L^{p(\cdot)} (\Omega)} \left\| \chi_{E_{n}} \right\|_{L^{\frac{p(\cdot) q(\cdot)}{p(\cdot) - q(\cdot)}} (\Omega)}. 
    \end{align*}

    Note that for each $n \in \mathbb{N}$ and each $\lambda > 0$,     
    \begin{displaymath}
        \chi_{\left\{ x \in \Omega : | \left( u \cdot \chi_{E_{n}} \right) (x) | > \lambda \right\}} = \chi_{\left\{ x \in \Omega : | u(x) | > \lambda \right\} \cap E_{n}} = \chi_{\left\{ x \in \Omega : | u(x) | > \lambda \right\} } \cdot \chi_{E_{n}}.
    \end{displaymath}

    Next, we estimate the following:     
    \begin{align*}
        & \left\| u \chi_{E_{n}} \right\|_{L^{q(\cdot), t} (\Omega)} = \left[ \int_{0}^{\infty} \left\| \chi_{\left\{ x \in \Omega : \left| \left( u \chi_{E_{n}} \right) (x) \right| > \lambda \right\}} \right\|_{L^{q(\cdot)} (\Omega)}^{t} \lambda^{t - 1} d \lambda \right]^{\frac{1}{t}} \\
        & = \left[ \int_{0}^{\infty} \left\| \chi_{\left\{ x \in \Omega : \left| u (x) \right| > \lambda \right\}} \cdot \chi_{E_{n}}\right\|_{L^{q(\cdot)} (\Omega)}^{t} \lambda^{t - 1} d \lambda \right]^{\frac{1}{t}} \\
        & \leq K \left[ \int_{0}^{\infty} \left\| \chi_{\left\{ x \in \Omega : \left| u (x) \right| > \lambda \right\}} \right\|_{L^{p(\cdot)} (\Omega)}^{t} \left\| \chi_{E_{n}} \right\|_{L^{\frac{p(\cdot) q(\cdot)}{p(\cdot) - q(\cdot)}} (\Omega)}^{t} \lambda^{t - 1} d \lambda \right]^{\frac{1}{t}} \\
        & = K \left\| \chi_{E_{n}} \right\|_{L^{\frac{p(\cdot) q(\cdot)}{p(\cdot) - q(\cdot)}} (\Omega)} \left[ \int_{0}^{\infty} \left\| \chi_{\left\{ x \in \Omega : \left| u (x) \right| > \lambda \right\}} \right\|_{L^{p(\cdot)} (\Omega)}^{t} \lambda^{t - 1} d \lambda \right]^{\frac{1}{t}} \\
        & = K \left\| \chi_{E_{n}} \right\|_{L^{\frac{p(\cdot) q(\cdot)}{p(\cdot) - q(\cdot)}} (\Omega)} \left\| u \right\|_{L^{p(\cdot), t} (\Omega)} \leq K \left\| \chi_{E_{n}} \right\|_{L^{\frac{p(\cdot) q(\cdot)}{p(\cdot) - q(\cdot)}} (\Omega)}.
    \end{align*}

    So, we obtain the following:
\begin{displaymath}
    \sup_{\left\| u \right\|_{L^{p(\cdot), t} (\Omega)} \leq 1} \left\| u \chi_{E_{n}} \right\|_{L^{q(\cdot), t} (\Omega)} \leq K \left\| \chi_{E_{n}} \right\|_{L^{\frac{p(\cdot) q(\cdot)}{p(\cdot) - q(\cdot)}} (\Omega)}.
\end{displaymath}

Since $| \Omega | < \infty$, $\chi_{\Omega}$ is integrable. By the dominated convergence theorem, $\left| E_{n} \right| = \int_{\Omega} \chi_{E_{n}} (x) d x \to 0$ as $n \to \infty$ because $\chi_{E_{n}} \to 0$ a.e as $n \to \infty$ and $\chi_{E_{n}} \leq \chi_{\Omega}$ for all $n \in \mathbb{N}$. Then, by assumption, 

\begin{displaymath}
    \lim_{n \to \infty} \sup_{\left\| u \right\|_{L^{p(\cdot), t} (\Omega)} \leq 1} \left\| u \chi_{E_{n}} \right\|_{L^{q(\cdot), t} (\Omega)} \leq K \lim_{n \to \infty} \left\| \chi_{E_{n}} \right\|_{L^{\frac{p(\cdot) q(\cdot)}{p(\cdot) - q(\cdot)}} (\Omega)} = 0.
\end{displaymath}
\end{proof}

To prepare for the lemma below, we adopt the following notation: If $a > 1$, then $a^{\infty} = \infty$.

\begin{lemma} \label{3.6}
     Let $p : \Omega \to [1, \infty)$, $q : \Omega \to [1, \infty)$ be measurable functions where $q(x) \leq p(x) \leq p_{+} < \infty$ for all $x \in \Omega$. Let $s : \Omega \to [1, \infty]$ be a function defined by:

\begin{displaymath}
s (x) :=
\begin{cases}
\frac{1}{p(x) - q(x)}
& \text{if } p(x) > q(x),\\
\infty & \text{if } p(x) = q(x). 
\end{cases} 
\end{displaymath}

     Assume that for all $a > 1$,

     \begin{equation} \label{finite}
         \int_{0}^{|\Omega|} a^{s^{*} (w)} \ d w < \infty.
     \end{equation}

     Then, for any sequence $(E_{n})_{n = 1}^{\infty}$ of measurable subsets of $\Omega$ such that $|E_{n}| \to 0$, we have

     \begin{displaymath}
         \left\| \chi_{E_{n}} \right\|_{L^{\frac{p(\cdot) q(\cdot)}{p(\cdot) - q(\cdot)}} (\Omega)} \to 0.
     \end{displaymath}
\end{lemma}
    
\begin{proof}
    By \eqref{finite}, it follows that $q (x) < p(x)$ for almost every $x \in \Omega$. Suppose that $\left( \left\| \chi_{E_{n}} \right\|_{L^{\frac{p(\cdot) q(\cdot)}{p(\cdot) - q(\cdot)}} (\Omega)} \right)_{n = 1}^{\infty}$ does not converge to $0$. Then, there exists $\epsilon > 0$ such that for all $k \in \mathbb{N}$, there exists $n_{k} \geq k$ where $\left\| \chi_{E_{n_{k}}} \right\|_{L^{\frac{p(\cdot) q(\cdot)}{p(\cdot) - q(\cdot)}} (\Omega)} \geq \epsilon$. Without loss of generality, we can assume that $\epsilon < 1$. By considering the subsequence $\left( \left\| \chi_{E_{n_{k}}} \right\|_{L^{\frac{p(\cdot) q(\cdot)}{p(\cdot) - q(\cdot)}} (\Omega)} \right)_{k = 1}^{\infty}$ of $\left( \left\| \chi_{E_{n}} \right\|_{L^{\frac{p(\cdot) q(\cdot)}{p(\cdot) - q(\cdot)}} (\Omega)} \right)_{n = 1}^{\infty}$, we see that for all $k \in \mathbb{N}$, $\left\| \chi_{E_{n_{k}}} \right\|_{L^{\frac{p(\cdot) q(\cdot)}{p(\cdot) - q(\cdot)}} (\Omega)} \geq \epsilon$. This means that for all $k \in \mathbb{N}$,

\begin{displaymath}
    \epsilon \leq \inf \left\{ \lambda > 0 : \int_{\Omega} \left| \frac{\chi_{E_{n_{k}}} (x)}{\lambda} \right|^{\frac{p(x) q(x)}{p(x) - q(x)}} d x \leq 1 \right\} = \inf \left\{ \lambda > 0 : \int_{\Omega} \left| \frac{\chi_{E_{n_{k}}} (x)}{\lambda} \right|^{p(x) q(x) s(x)} d x \leq 1 \right\}.
\end{displaymath}
    
Let $\delta \in (0, \epsilon)$. By the definition of infimum, we deduce that

\begin{displaymath}
    \int_{E_{n_{k}}} \left( \frac{1}{\delta} \right)^{p(x) q(x) s(x)} d x = \int_{\Omega} \left| \frac{\chi_{E_{n_{k}}} (x)}{\delta} \right|^{p(x) q(x) s(x)} d x > 1.
\end{displaymath}

Next, we apply the Hardy-Littlewood inequality to obtain the following estimate:

\begin{align*}
    & \int_{E_{n_{k}}} \left( \frac{1}{\delta} \right)^{p(x) q(x) s(x)} d x \leq \int_{E_{n_{k}}} \left( \frac{1}{\delta^{p_{+} q_{+}}} \right)^{s(x)} d x \text{ (since } \delta < 1) \\
    & = \int_{\Omega} \left( \frac{1}{\delta^{p_{+} q_{+}}} \right)^{s(x)} \chi_{E_{n_{k}}} (x) d x \\
    & \leq \int_{0}^{\infty} \left[ \left( \frac{1}{\delta^{p_{+} q_{+}}} \right)^{s(\cdot)} \right]^{*} (y) \left( \chi_{E_{n_{k}}} \right)^{*} (y) d y \text{ (by the Hardy-Littlewood inequality)} \\
    & = \int_{0}^{\infty} \left( \frac{1}{\delta^{p_{+} q_{+}}} \right)^{s^{*} (y)} \chi_{ \left[0, \left| E_{n_{k}} \right| \right)} (y)  \ d y \text{ (by \thref{exponent})} \\
    & = \int_{0}^{\left| E_{n_{k}}\right|} \left( \frac{1}{\delta^{p_{+} q_{+}}} \right)^{s^{*} (y)} \ d y.
\end{align*}

    Combining the above two inequalities, we conclude that for all $k \in \mathbb{N}$, 
    \begin{align} \label{contradiction}
        \int_{0}^{\left| E_{n_{k}}\right|} \left( \frac{1}{\delta^{p_{+} q_{+}}} \right)^{s^{*} (y)} \ d y > 1 .
    \end{align}
    
    Meanwhile, since $\frac{1}{\delta^{p_{+} q_{+}}} > 1$ by assumption, we know that 

    \begin{displaymath}
        \int_{0}^{\left| \Omega \right|} \left( \frac{1}{\delta^{p_{+} q_{+}}} \right)^{s^{*} (y)} \ d y < \infty.
    \end{displaymath}

    So, by the absolute continuity of the Lebesgue integral, there exists $\beta > 0$ such that for all $\omega < \beta$, 

    \begin{displaymath}
        \int_{0}^{\omega} \left( \frac{1}{\delta^{p_{+} q_{+}}} \right)^{s^{*} (y)} \ d y < 1.
    \end{displaymath}
    By assumption, $\left| E_{n_{k}} \right| \to 0$. So, there exists $K \in \mathbb{N}$ such that for all $k \geq K$, $\left| E_{n_{k}} \right| < \beta$. This implies that for all $k \geq K$,

    \begin{displaymath}
        \int_{0}^{\left| E_{n_{k}}\right|} \left( \frac{1}{\delta^{p_{+} q_{+}}} \right)^{s^{*} (y)} \ d y < 1.
    \end{displaymath}

    This contradicts \eqref{contradiction}.
\end{proof}

Now we can state a theorem which is an extension of \cite[Theorem 3.4]{EGN}.

\begin{theorem} \thlabel{almostcompact}
     Let $p : \Omega \to [1, \infty)$, $q : \Omega \to [1, \infty)$ be measurable functions where $q(x) \leq p(x) \leq p_{+} < \infty$ for all $x \in \Omega$. Let $s : \Omega \to [1, \infty]$ be a function defined by:

\begin{displaymath}
s (x) :=
\begin{cases}
\frac{1}{p(x) - q(x)}
& \text{if } p(x) > q(x),\\
\infty & \text{if } p(x) = q(x). 
\end{cases} 
\end{displaymath}

     Assume that for all $a > 1$,

     \begin{displaymath}
         \int_{0}^{|\Omega|} a^{s^{*} (w)} \ d w < \infty.
     \end{displaymath}

     Then, $L^{p(\cdot), t} (\Omega)$ is almost compactly embedded into $L^{q(\cdot), t} (\Omega)$ for any  $t \in [1, \infty)$.
\end{theorem}
\begin{proof}
This follows from Lemma \ref{3.5} and Lemma \ref{3.6}.    
\end{proof}

\begin{lemma} \label{3.8}
    Let $t\in [1,\infty)$ and $1 \leq p < d$, where $d \geq 2$ is a positive integer. Let $\Omega$ be a bounded domain on $\mathbb{R}^{d}$. Let $q : \Omega \to [1, \infty)$ be a measurable function satisfying $1 \leq q(x) \leq \frac{d p}{d - p}$ for almost every $x \in \Omega$. Suppose that there exist $x_{0} \in \Omega$, $C > 0$, $\eta > 0$ and $0 < \ell < 1$ such that 

    \begin{displaymath}
        \text{ess } \sup \left\{ q(x) \in [1, \infty) : x \in \Omega \text{ and } \left| x - x_{0} \right| \geq \eta \right\} < \frac{d p}{d - p}
    \end{displaymath}

    and

    \begin{displaymath}
        q(x) \leq \frac{d p}{d - p} - \frac{C}{\left| \log \left( \frac{1}{\left| x - x_{0} \right|} \right) \right|^{\ell}} \text{ for a.e. } x \in \Omega \text{ with } \left| x - x_{0} \right| \leq \eta.
    \end{displaymath}

    Then, $L^{p^{*}, t} (\Omega)$ is almost compactly embedded into $L^{q(\cdot), t} (\Omega)$.
\end{lemma}

\begin{proof}
    Without loss of generality, we can assume that $|x_{0}| = 0 \in \Omega$. Also, we can assume that there exist  $C > 0$, $1>\eta > 0$ and $0 < \ell < 1$ such that $q(x) \leq r(x)$, where $r : \Omega \to [1, \infty]$ is defined by:

\begin{displaymath}
r (x) :=
\begin{cases}
\frac{d p}{d - p} - \frac{C}{\left| \log \left( \frac{1}{\left| x \right|} \right) \right|^{\ell}}
& \text{if }|x| \leq \eta,\\
\frac{d p}{d - p} - \frac{C}{\left| \log \left( \frac{1}{\eta} \right) \right|^{\ell}} & \text{if } |x| > \eta. 
\end{cases} 
\end{displaymath}

So, $\frac{1}{p^{*} - q(x)} \leq \frac{1}{p^{*} - r(x)}$ for all $x \in \Omega$. It follows that $\left[ \frac{1}{p^{*} - q(\cdot)} \right]^{*} (t) \leq \left[ \frac{1}{p^{*} - r(\cdot)} \right]^{*} (t)$ for all $t > 0$. Thus, it remains to show that for all $a > 1$,

\begin{displaymath}
    \int_{0}^{| \Omega |} a^{\left[ \frac{1}{p^{*} - r(\cdot)} \right]^{*} (t)}  \ dt < \infty. 
\end{displaymath}

For the simplicity of notation, we let $s(x) := \frac{1}{p^{*} - r(x)}$. So,

\begin{displaymath}
s (x) :=
\begin{cases}
\frac{1}{C}\left| \log \left( \frac{1}{\left| x \right|} \right) \right|^{\ell}
& \text{if }|x| \leq \eta,\\
\frac{1}{C}\left| \log \left( \frac{1}{\eta} \right) \right|^{\ell} & \text{if } |x| > \eta. 
\end{cases} 
\end{displaymath}

Recalling the definition for the distribution function, for each $\alpha > 0$,   we define
$$
d_{s} (\alpha) := \left| \left\{ x \in \Omega : | s(x) | > \alpha \right\} \right|.$$ 

For $0 \leq \alpha < \frac{1}{C}\left| \log \left( \frac{1}{\eta} \right) \right|^{\ell}$ we have obviously  $d_{s} (\alpha) = | \Omega |$. 

Meanwhile, for $\alpha \geq \frac{1}{C}\left| \log \left( \frac{1}{\eta} \right) \right|^{\ell}$, we get
\begin{displaymath}
    d_{s} (\alpha) = \left| \left\{ x \in \Omega :  \frac{1}{C}\left| \log \left( \frac{1}{\left| x \right|} \right) \right|^{\ell}  > \alpha \right\} \right| = \left| \left\{ x \in \Omega : \left| \log \left( \frac{1}{\left| x \right|} \right) \right|  > \left( C \alpha \right)^{\frac{1}{\ell}} \right\} \right|
\end{displaymath}
\begin{displaymath}
    = \left| \left\{ x \in \Omega : - \log \left( \left| x \right| \right)  > \left( C \alpha \right)^{\frac{1}{\ell}} \right\} \right| = \left| \left\{ x \in \Omega : |x| < e^{- \left( C \alpha \right)^{\frac{1}{\ell}} } \right\} \right| = \nu_{d} e^{- d \left( C \alpha \right)^{\frac{1}{\ell}} }. 
\end{displaymath}

Summarizing the above, we obtain:

\begin{displaymath}
d_{s} (\alpha) :=
\begin{cases}
| \Omega | & \text{if } 0 \leq \alpha < \frac{1}{C}\left| \log \left( \frac{1}{\eta} \right) \right|^{\ell},\\
\nu_{d} e^{- d \left( C \alpha \right)^{\frac{1}{\ell}} }  & \text{if } \alpha \geq \frac{1}{C}\left| \log \left( \frac{1}{\eta} \right) \right|^{\ell}. 
\end{cases} 
\end{displaymath}

Next, recall the definition of the non-increasing rearrangement function. For $t \geq | \Omega |$, we have $s^{*} (t) = 0$, and for $\nu_{d} \eta^{d} \leq t < | \Omega|$, $s^{*} (t) = \frac{1}{C} \left| \log \left( \frac{1}{\eta} \right)^{\ell} \right|$. Finally, let $0 \leq t < \nu_{d} \eta^{d}$, then,
\begin{displaymath}
    d_{s} (\beta) \leq t \Longleftrightarrow \nu_{d} e^{- d \left( C \beta \right)^{\frac{1}{\ell}}} \leq t \Longleftrightarrow \beta \geq \frac{1}{C d^{\ell}} \left[ \log \left( \frac{\nu_{d}}{t} \right) \right]^{\ell}.
\end{displaymath}

So, we deduce that 
\begin{displaymath}
s^{*} (t) :=
\begin{cases}
\frac{1}{C d^{\ell}} \left[ \log \left( \frac{\nu_{d}}{t} \right) \right]^{\ell} & \text{if } 0 \leq t < \nu_{d} \eta^{d},\\
\frac{1}{C} \left| \log \left( \frac{1}{\eta} \right) \right|^{\ell}  & \text{if } \nu_{d} \eta^{d} \leq t < | \Omega |, \\
0 & \text{ if } t \geq | \Omega |.
\end{cases} 
\end{displaymath}

Finally, we proceed to show that for each $a > 1$, $\int_{0}^{\nu_{d} \eta^{d}} a^{s^{*} (t)} \ d t < \infty$.

\begin{displaymath}
    \int_{0}^{\nu_{d} \eta^{d}} a^{\frac{1}{C d^{\ell} } \left[ \log \left( \frac{\nu_{d}}{t} \right) \right]^{\ell}} d t = \int_{0}^{\nu_{d} \eta^{d}} e^{\frac{1}{C d^{\ell}} \left[ \log \left( \nu_{d}  \right) - \log (t) \right]^{\ell} \log(a)} d t
\end{displaymath}
    
\begin{displaymath}
    = \int_{\infty}^{- d \log (\eta) } - e^{\frac{1}{C d^{\ell}} y^{\ell} \log(a)} e^{\log \left( \nu_{d} \right) - y} d y = \nu_{d} \int_{- d \log (\eta) }^{\infty} e^{\frac{\log(a)}{C d^{\ell}} y^{\ell} - y} d y. 
\end{displaymath}

Since $0 < \ell < 1$, there exists $0 < \omega < 1$ such that for all $y \geq - d \log (\eta)$, $e^{\frac{\log(a)}{C d^{\ell}} y^{\ell} - y} \leq e^{- \omega y}$. Hence,

\begin{displaymath}
    \int_{0}^{\nu_{d} \eta^{d}} a^{\frac{1}{C d^{\ell} } \left[ \log \left( \frac{\nu_{d}}{t} \right) \right]^{\ell}} d t \leq \nu_{d} \int_{- d \log (\eta) }^{\infty} e^{- \omega y} d y < \infty.
\end{displaymath}

Thus, for all $a > 1$

\begin{displaymath}
    \int_{0}^{| \Omega |} a^{\left[ \frac{1}{p^{*} - q (\cdot) } \right]^{*} (t)} \ d t \leq \int_{0}^{| \Omega |} a^{\left[ \frac{1}{p^{*} - r (\cdot) } \right]^{*} (t)} \ d t < \infty.
\end{displaymath}

By \thref{almostcompact}, $L^{p^{*}, t} (\Omega)$ is almost compactly embedded into $L^{q(\cdot), t} (\Omega)$.
\end{proof}

Now we are ready to state the conditions on $q(\cdot)$ which guarantees compactness. 

\begin{theorem} \thlabel{q-compact}
    Let $1 \leq p < d$, where $d \geq 2$ is a positive integer. Let $\Omega$ be a bounded domain on $\mathbb{R}^{d}$. Let $q : \Omega \to [1, \infty)$ be a measurable function satisfying $1 \leq q(x) \leq \frac{d p}{d - p}$ for almost every $x \in \Omega$. Suppose that there exist $x_{0} \in \Omega$, $C > 0$, $\eta > 0$ and $0 < \ell < 1$ such that 
    \begin{displaymath}
        \text{ess } \sup \left\{ q(x) \in [1, \infty) : x \in \Omega \text{ and } \left| x - x_{0} \right| \geq \eta \right\} < \frac{d p}{d - p}
    \end{displaymath}
    and
    \begin{displaymath}
        q(x) \leq \frac{d p}{d - p} - \frac{C}{\left| \log \left( \frac{1}{\left| x - x_{0} \right|} \right) \right|^{\ell}} \text{ for a.e. } x \in \Omega \text{ with } \left| x - x_{0} \right| \leq \eta.
    \end{displaymath}

    Then, the embedding $I : W_0^{1, p} (\Omega) \to L^{q(\cdot), p} (\Omega)$ is compact.
\end{theorem}
\begin{proof}
    It follows instantly from Proposition \ref{prop3.4} and Lemma \ref{3.8}.
\end{proof}

The next statement provides a condition on $q(\cdot)$ under which Sobolev embedding into $L^{q(\cdot),p}$ is non-compact. 
The idea behind the proof is inspired by [\cite{KS2008}, Theorem 1], and from the above, we can see that the condition on $q(\cdot)$ appears to be quite sharp (in light of the logarithmic term's power).

\begin{theorem} \thlabel{propnon-compact}
	Let $1 \leq p < d$, where $d \geq 2$ is a positive integer, and $\Omega$ be a bounded domain on $\mathbb{R}^{d}$. Let $q : \Omega \to [1, \infty)$ be a measurable function satisfying $1 \leq q(x) \leq \frac{d p}{d - p}$ for almost every $x \in \Omega$. Suppose that there exist $x_{0} \in \Omega$ and constants $C_{0}$, $\eta_{0} > 0$ such that 
	
	\begin{equation} \label{non-compact}
		q(x) \geq \frac{d p}{d - p} - \frac{C_{0}}{ \left| \log \left( \frac{1}{| x - x_{0} |} \right) \right|} \text{ for almost every } x \in \Omega \text{ with } | x - x_{0} | < \eta_{0}.
	\end{equation}
	
	Then, the embedding $I : W_0^{1, p} (\Omega) \to L^{q(\cdot), p} (\Omega)$ is not compact. 
\end{theorem}

\begin{proof}
	Without loss of generality, we can assume $x_{0} = 0$. Define $r : \Omega \to [0, \infty)$ by $r(x) := \frac{d p}{d - p} - q(x)$. Now, suppose that $E$ is compact. Define the bump function $\phi \in C_{0}^{\infty} (\mathbb{R}^{d})$ such that $\phi(x) = 1$ for all $| x | \leq \frac{1}{2}$ and $\phi(x) = 0$ for all $| x | > 1$. 
	
	\bigskip
	
	For each $n \in \mathbb{N}$, define $\phi_{n} : \Omega \to [0, \infty)$ by $\phi_{n} (x) := n^{\frac{d - p}{p}} \phi( n x )$. We observe that for sufficiently large $n \in \mathbb{N}$, $\phi_{n} \in C_{0}^{\infty} (\Omega)$. Also, for each $n \in \mathbb{N}$, we have the following two equalities:
		\begin{align*}
		& \int_{\Omega} | ( \nabla \phi_{n} ) (x) |^{p} \ d x = \int_{B_{\frac{1}{n}} (0)} \left[ \sum_{j = 1}^{d} \left| \frac{\partial \phi_{n}}{\partial x_{j}} (x) \right|^{2} \right]^{\frac{p}{2}} d x = \int_{B_{{1}} (0)} \left[ \sum_{j = 1}^{d} \left| n^{\frac{d - p}{p} + 1}  \frac{\partial \phi}{\partial x_{j}} (y) \right|^{2} \right]^{\frac{p}{2}} {d y \over n^d} \\
		& = \frac{n^{d}}{n^{d}} \int_{B_{1} (0)} \left[ \sum_{j = 1}^{d} \left| \frac{\partial \phi}{\partial x_{j}} (y) \right|^{2} \right]^{\frac{p}{2}} d y = \int_{B_{1} (0)} | ( \nabla \phi ) (y) |^{p} \ d y,
	\end{align*}	
	and
	\begin{equation} 	\label{L_p phi}
		\int_{\Omega} | \phi_{n} (x) |^{p} \ d x = \int_{B_{\frac{1}{n}} (0)} n^{d - p} | \phi (n x) |^{p} \ d x = \frac{1}{n^{p}} \int_{B_{1} (0)} | \phi (y) |^{p} \ d y.
	\end{equation}
	
	Now, let us estimate $\left\| \chi_{\left\{ z \in \Omega : \phi(n z) = 1 \right\} }  \right\|_{L^{q (\cdot)}}$. Choose $1>\delta > 0$ sufficiently small such that $ \delta < e^{- C_{0} \frac{d - p}{p}} {\nu_{d}}/{2^{d}},$ where $\nu_{d} =2^d \left| B_{\frac{1}{2}}(0) \right|$.

	\begin{align*}
		& \int_{\Omega} \left| \frac{n^{\frac{d - p}{p}}}{\delta} \chi_{\left\{ z \in \Omega : \phi(n z) = 1 \right\} } (x) \right|^{q(x)} \ d x \geq \frac{1}{\delta} \int_{\Omega} \left| n^{\frac{d - p}{p}} \chi_{\left\{ z \in \Omega : \phi(n z) = 1 \right\} } (x) \right|^{q(x)} \ d x \\
		& = \frac{1}{\delta} \int_{B_{\frac{1}{2 n}} (0)} n^{\frac{d - p}{p} \left[ \frac{d p}{d - p} - r(x) \right] } \ d x  = \frac{n^{d}}{\delta} \int_{B_{\frac{1}{2 n}} (0)} n^{- \frac{d - p}{p} r(x) } \ d x \\
		& = \frac{1}{\delta} \int_{B_{\frac{1}{2}} (0)} n^{- \frac{d - p}{p} r(\frac{y}{n}) } \ dy  \geq \frac{1}{\delta} \int_{B_{\frac{1}{2}} (0)} \frac{1}{n^{\frac{C_{0} (d - p)}{p \left| \log \left( \frac{|y|}{n} \right) \right|}}} \ dy \\
		& = \frac{1}{\delta} \int_{B_{\frac{1}{2}} (0)} e^{- \frac{d - p}{p} \frac{C_{0} \log(n) }{\left| \log \left( \frac{|y|}{n} \right) \right|}} \ dy  = \frac{1}{\delta} \int_{B_{\frac{1}{2}} (0)} e^{- \frac{d - p}{p} \frac{C_{0} \log(n) }{\left| \log (n) - \log (|y|) \right|}} \ d y \\
		& \geq \frac{1}{\delta} \int_{B_{\frac{1}{2}}(0)} e^{- \frac{d - p}{p} \frac{C_{0} \log(n) }{ \log (n) }} \ d y = e^{- C_{0} \frac{d - p}{p}} \frac{\nu_{d}}{ 2^{d} \delta} > 1.
	\end{align*}
	
	
	

\begin{displaymath}
    \text{Since }\rho_{q(\cdot)} \left[ \frac{ \chi_{\left\{ z \in \Omega : \phi (n z) = 1\right\} } }{ \delta n^{- \frac{d - p}{p}}}  \right] > 1, \quad \left\| \chi_{\left\{ z \in \Omega : \phi(n z) = 1 \right\} }  \right\|_{L^{q (\cdot)}} \geq \frac{\delta}{n^{\frac{d - p}{p}}}.
\end{displaymath}

	Next, for each $n \in \mathbb{N}$,
	\begin{align*}
		& \left\| \phi_{n} \right\|_{L^{q(\cdot), p}}^{p} = \int_{0}^{\infty} \lambda^{p - 1} \left\| \chi_{\{ z \in \Omega : | \phi_{n} (z) | > \lambda \} } \right\|_{L^{q(\cdot)}}^{p} d \lambda = \int_{0}^{\infty} \lambda^{p - 1} \left\| \chi_{ \left\{ z \in \Omega : n^{\frac{d - p}{p}} \phi( n z ) > \lambda \right\} } \right\|_{L^{q(\cdot)}}^{p} d \lambda \\
		& = \int_{0}^{\infty} \lambda^{p - 1} \left\| \chi_{ \left\{ z \in \Omega : | \phi (n z) | > \frac{\lambda}{n^{\frac{d - p}{p}}} \right\} } \right\|_{L^{q(\cdot)}}^{p} \ d \lambda  = n^{d - p} \int_{0}^{\infty} \tilde{\lambda}^{p - 1} \left\| \chi_{\left\{ z \in \Omega : \phi(n z) > \tilde{\lambda} \right\} }  \right\|_{L^{q (\cdot)}}^{p} \ d \tilde{\lambda} \\
		& = n^{d - p} \int_{0}^{1} \tilde{\lambda}^{p - 1} \left\| \chi_{\left\{ z \in \Omega : \phi(n z) > \tilde{\lambda} \right\} }  \right\|_{L^{q (\cdot)}}^{p} \ d \tilde{\lambda}  \geq n^{d - p} \int_{0}^{1} \tilde{\lambda}^{p - 1} 
		\left\| \chi_{\left\{ z \in \Omega : \phi(n z) = 1 \right\} }  \right\|_{L^{q (\cdot)}}^{p} \ d \tilde{\lambda} \\
		& \geq n^{d - p} \int_{0}^{1} \tilde{\lambda}^{p - 1} \frac{\delta^{p}}{n^{d - p}} 
		\ d \tilde{\lambda} = \frac{\delta^{p}}{p} > 0.
	\end{align*}
	
	
	
	
	Since $\left\| \nabla \phi_{n} \right\|_p \leq C$ for all $n \in \mathbb{N}$, by the assumption that $E$ is compact, there exist a subsequence $(\phi_{n_{k}})_{k = 1}^{\infty}$ of $(\phi_{n})_{n = 1}^{\infty}$ and $\psi \in L^{q(\cdot), p} (\Omega)$ such that $\left\| \phi_{n_{k}} - \psi \right\|_{L^{q(\cdot), p}} \to 0$. In particular, there exists another subsequence $(\phi_{n_{k_{j}}})_{j = 1}^{\infty}$ of $(\phi_{n_{k}})_{k = 1}^{\infty}$ that converges to $\psi$ almost everywhere. By the above, for all $n \in \mathbb{N}$, $\left\| \phi_{n} \right\|_{L^{q(\cdot), p}} \geq \frac{\delta}{p^{\frac{1}{p}}}$. So, $\psi \not = 0$.
	Meanwhile, by (\ref{L_p phi})  we obtain $\left\| \phi_{n} \right\|_{L^{p}(\Omega) } \to 0$ and this contradicts the fact that $\psi \not = 0$.
	
	
\end{proof}

We note that under condition (\ref{non-compact}), the Sobolev embedding is  only non-compact at the neighborhood of the point $x_0 \in \Omega$ and the restriction of the Sobolev embedding at any region of $\Omega$ not containing a neighborhood of $x_0$ is compact.

\section{Quality of Non-Compactness}

In this section, we show that condition (\ref{non-compact}) will give us a Sobolev embedding which is almost as non-compact as the Sobolev embedding in (\ref{2}), i.e. non-compactness concentrated at just one point of the domain could produce non-compactness comparable to the ``most'' non-compact Sobolev embedding  (\ref{2}) which is non-compact at the neighborhood of each point of the domain.

Let us recall that  in  \thref{propnon-compact},
the embedding $I : W^{1, p} (\Omega) \to L^{q(.), p} (\Omega)$ is non-compact since the function $q(\cdot)$ approaches the value $p^{*}$ at a fast enough rate at $x_0\in \Omega$. So, it makes sense to introduce the following quantities:

\begin{equation} \label{gamma}
	\gamma_{r} := \sup_{f \in W_0^{1, p} [\Omega \cap B_{r} (x_0) ]} \frac{\| f \|_{L^{q(\cdot), p} } }{\| \nabla f \|_{L^{p}} } \ge   \lim_{r \to 0^{+}} \sup_{f \in W_{0}^{1, p} [\Omega \cap B_{r} (x_0) ]} \frac{\| f \|_{L^{q(\cdot), p} } }{\| \nabla f \|_{L^{p}}  }=:\gamma.
\end{equation}


\begin{lemma} \thlabel{gamma lemma}
	 Suppose that the conditions of  \thref{propnon-compact} are satisfied. Then,  $\gamma_r \ge  \gamma > 0$. 
\end{lemma}

\begin{proof}
Without a loss of generality, we can assume, in the rest of this proof, that we have $\Omega = B_1(0)$ and $x_0=0$.	

As in the proof of \thref{propnon-compact}, we first consider a bump function $ \phi \in C_{0}^{\infty} (\mathbb{R}^{d})$ such that $\phi(x) = 1$ for all $| x | \leq \frac{1}{2}$ and $\phi(x) = 0$ for all $| x | > 1$.

	For each $r > 0$, define $\phi_{r} : \Omega \to [0, \infty)$ by $\phi_{r} (x) := \left( \frac{1}{r} \right)^{\frac{d - p}{p}} \phi \left( \frac{1}{r} x \right)$. Following the calculations in the proof of \thref{propnon-compact}, we observe that for all $r > 0$,	
	\begin{displaymath}
		\left\| \nabla \phi_{r} \right\|_{L^{p} [B_{r} (0)]} = \left[ \int_{B_{r} (0)} \left| \left( \nabla \phi_{r} \right) (x) \right|^{p} d x \right]^{\frac{1}{p}} = \left[ \int_{B_{1} (0)} \left| \left( \nabla \phi \right) (x) \right|^{p} (x) d x \right]^{\frac{1}{p}} = \left\| \nabla \phi \right\|_{L^{p} [B_{1} (0)]}
	\end{displaymath}
	
	Similarly, there exists $0 < \delta < 1$ such that for all $r \in (0, 1)$,	
	\begin{displaymath}
		\left\| \phi_{r} \right\|_{L^{q (\cdot), p} (\Omega)} \geq \frac{\delta}{p^{\frac{1}{p}}} > 0.
	\end{displaymath}

	Thus, by definition, for all $r \in (0, 1)$,	
	\begin{displaymath}
		\gamma_{r}:= \sup_{f \in W^{1, p} [B_{r} (0) ]} \frac{\| f \|_{L^{q(\cdot), p} [B_{r} (0)]} }{\| \nabla f \|_{L^{p}[B_{r} (0)]} } \geq \frac{\delta}{p^{\frac{1}{p}} \left\| \nabla \phi \right\|_{L^{p} [B_{1} (0)]}} > 0.
	\end{displaymath}
	
	\begin{displaymath}
		\text{Therefore, } \gamma := \lim_{r \to 0^{+}} \gamma_{r} \geq \frac{\delta}{p^{\frac{1}{p}} \left\| \nabla \phi \right\|_{L^{p} [B_{1} (0)]}} > 0.
	\end{displaymath}
\end{proof}


Now, we introduce a couple of simple technical lemmas.

\begin{lemma} \thlabel{symmetrization}
Let $\Omega \subseteq \mathbb{R}^{d}$ be an open set such that $| \Omega | \leq 1$, $1 \leq p < \infty$, and  $q : \Omega \to [1, \infty)$. Suppose that there exists a decreasing function $\tilde{q} : [0, \infty) \to [1, \infty)$ and $x_0\in \Omega$ such that $q (x) = \tilde{q} \left( \left| x - x_{0} \right| \right)$ for all $x \in \Omega$.  $q : \Omega \to [1, \infty).$ 

Then, given any function $f : \Omega \to \mathbb{C}$, we have that 

\begin{displaymath}
    \left\| f \right\|_{L^{q(\cdot), p} (\Omega)} \leq \left\| f^{\#} \right\|_{L^{q(\cdot), p} (\mathbb{R}^{d})}.
\end{displaymath}
\end{lemma}

\begin{proof}
    Without loss of generality, we can assume that $x_{0} = 0$. First, we observe that for each $\lambda > 0$,     
    \begin{displaymath}
        \left\| \chi_{\left\{ x \in \mathbb{R}^{d} : f^{\#} (x) > \lambda \right\} } \right\|_{L^{q(\cdot)} (\mathbb{R}^{d})} \leq \left\| \chi_{\Omega^{\#}} \right\|_{L^{q(\cdot)} (\mathbb{R}^{d})} \leq 1 \text{ since } \int_{\mathbb{R}^{d}} \left| \chi_{\Omega^{\#}} (x) \right|^{q(x)} d x = \left| \Omega^{\#} \right| = \left| \Omega \right| \leq 1.
    \end{displaymath}

    Let $\lambda > 0$. Note that there exists $r_{\lambda} > 0$ such that $\chi_{\left\{ x \in \mathbb{R}^{d} : f^{\#} (x) > \lambda \right\} } = \chi_{B_{r_{\lambda}} (0) }$.

    By the equimeasurability of the symmetric decreasing rearrangement, we have that:

    \begin{displaymath} 
        \left| \left\{ x \in \mathbb{R}^{d} : f^{\#} (x) > \lambda \right\} \right| = \left| \left\{ x \in \Omega : | f(x) | > \lambda \right\} \right|.
    \end{displaymath}

Denote 
\begin{align*}
    & A_1:=  \left\{ x \in \mathbb{R}^{d} : f^{\#} (x) > \lambda, |f(x)| > \lambda \right\} , \qquad
    A_2:= \left\{ x \in \mathbb{R}^{d} : f^{\#} (x) > \lambda, |f(x)| \leq \lambda \right\} \\
    & A_3 :=  \left\{ x \in \Omega : | f(x) | > \lambda, f^{\#} (x) \leq \lambda \right\} .
\end{align*}

\begin{equation} \label{same measure}
    \text{Since } \left| A_{1} \cup A_{2} \right| = \left| A_{1} \cup A_{3} \right|, A_{1} \cap A_{2} = \emptyset \text{ and } A_{1} \cap A_{3} = \emptyset, \text{ we see that } \left| A_{2} \right| = \left| A_{3}. \right|.
\end{equation}

    Now, we show that $\left\| \chi_{\left\{ x \in \Omega : |f(x)| > \lambda \right\} } \right\|_{L^{q(\cdot)} (\Omega)} \leq \left\| \chi_{\left\{ x \in \Omega : f^{\#} (x) > \lambda \right\} } \right\|_{L^{q(\cdot)} (\mathbb{R}^{d})}$.
    
    \begin{align*}
        & \int_{\Omega} \left| \frac{ \chi_{\left\{ x \in \Omega : |f(x)| > \lambda \right\} } }{\left\| \chi_{ \left\{ y \in \mathbb{R}^{d}: f^{\#} (y) > \lambda \right\} }\right\|}_{q(\cdot)} \right|^{q(x)} d x = \int_{\left\{ x \in \Omega : |f(x)| > \lambda \right\}} \frac{ 1 }{\left\| \chi_{ \left\{ y \in \mathbb{R}^{d}: f^{\#} (y) > \lambda \right\} } \right\|_{q(\cdot)}^{q(x)} } d x \\
        & = \int_{A_{1}} \frac{ 1 }{\left\| \chi_{ \left\{ y \in \mathbb{R}^{d}: f^{\#} (y) > \lambda \right\} }\right\|_{q(\cdot)}^{q(x)} } d x + \int_{A_{3}} \frac{ 1 }{\left\| \chi_{ \left\{ y \in \mathbb{R}^{d}: f^{\#} (y) > \lambda \right\} } \right\|_{q(\cdot)}^{q(x)} } d x \\
        & = \int_{A_{1}} \frac{ 1 }{\left\| \chi_{ \left\{ y \in \mathbb{R}^{d}: f^{\#} (y) > \lambda \right\} }\right\|_{q(\cdot)}^{q(x)} } d x + \int_{\left\{ x \in \mathbb{R}^{d} \setminus B_{r_{\lambda}} (0) : |f(x)| > \lambda, f^{\#} (x) \leq \lambda \right\}} \frac{ 1 }{\left\| \chi_{ \left\{ y \in \mathbb{R}^{d}: f^{\#} (y) > \lambda \right\} } \right\|_{q(\cdot)}^{q(x)} } d x \\
        & \leq \int_{A_{1}} \frac{ 1 }{\left\| \chi_{ \left\{ y \in \mathbb{R}^{d}: f^{\#} (y) > \lambda \right\} }\right\|_{q(\cdot)}^{q(x)} } d x + \int_{\left\{ x \in \mathbb{R}^{d} \setminus B_{r_{\lambda} } (0) : |f(x)| > \lambda, f^{\#} (x) \leq \lambda \right\}} \frac{ 1 }{\left\| \chi_{ \left\{ y \in \mathbb{R}^{d}: f^{\#} (y) > \lambda \right\} } \right\|_{q(\cdot)}^{\tilde{q} (r_{\lambda}) } } d x \\
        & = \int_{A_{1}} \frac{ 1 }{\left\| \chi_{ \left\{ y \in \mathbb{R}^{d}: f^{\#} (y) > \lambda \right\} }\right\|_{q(\cdot)}^{q(x)} } d x + \int_{\left\{ x \in \mathbb{R}^{d} : f^{\#} (x) > \lambda, | f(x) | \leq \lambda \right\}} \frac{ 1 }{\left\| \chi_{ \left\{ y \in \mathbb{R}^{d}: f^{\#} (y) > \lambda \right\} } \right\|_{q(\cdot)}^{\tilde{q} (r_{\lambda}) } } d x \text{ by \eqref{same measure}} \\
        & \leq \int_{A_{1}} \frac{ 1 }{\left\| \chi_{ \left\{ y \in \mathbb{R}^{d}: f^{\#} (y) > \lambda \right\} }\right\|_{q(\cdot)}^{q(x)} } d x + \int_{A_{2}} \frac{ 1 }{\left\| \chi_{ \left\{ y \in \mathbb{R}^{d}: f^{\#} (y) > \lambda \right\} } \right\|_{q(\cdot)}^{q(x)} } d x \\
        & = \int_{\left\{ x \in \mathbb{R}^{d} : f^{\#} (x) > \lambda \right\}} \frac{ 1 }{\left\| \chi_{ \left\{ y \in \mathbb{R}^{d}: f^{\#} (y) > \lambda \right\} }\right\|_{q(\cdot)}^{q(x)} } d x = 1.
    \end{align*}      

    Finally, we see that
    \begin{align*}
        & \left\| f \right\|_{L^{q(\cdot), p} (\Omega)} = \left[ \int_{0}^{\infty} \lambda^{p - 1} \left\| \chi_{\left\{ x \in \Omega : |f(x)| > \lambda \right\} } \right\|_{L^{q(\cdot)}}^{p} d \lambda \right]^{\frac{1}{p}} \\
        & \leq \left[ \int_{0}^{\infty} \lambda^{p - 1} \left\| \chi_{\left\{ x \in \mathbb{R}^{d} : f^{\#} (x) > \lambda \right\} } \right\|_{L^{q(\cdot)}}^{p} d \lambda \right]^{\frac{1}{p}} \leq \left\| f^{\#} \right\|_{L^{q(\cdot), p} (\mathbb{R}^{d})}.
    \end{align*}
\end{proof}

\begin{lemma} \thlabel{main}
	Let $p \in [1, \infty)$, $\varepsilon > 0$. Then, there exists $k_{\varepsilon} \in \mathbb{N}$ such that for all sequences $(\alpha_{n})_{n = 1}^{\infty} \subseteq \mathbb{C}$,
	
	\begin{equation}
		\sum_{n \in \left\{ m \in \mathbb{N} : | \alpha_{m}| \leq \frac{1}{2^{m + k_{\varepsilon}}} \right\} } | \alpha_{n} |^{p} < \varepsilon
		\sum_{n=1}^\infty |\alpha_n |^p.
		\label{series}
	\end{equation}
	
\end{lemma}
\begin{proof}
	Setting $k_{\varepsilon}$ such that 
	\[  \sum_{n=1}^{\infty}  \frac{1}{2^{m + k_{\varepsilon}}} < \varepsilon
	\]
	gives us (\ref{series}) instantly.
\end{proof}

We skip the proof of the next lemma as it is obvious.
\begin{lemma} \thlabel{almost orthogonal}
	Let $X$ be a Banach space. Let $p \in [1, \infty]$, $\varepsilon > 0$. Suppose that $(\alpha_{n})_{n = 1}^{\infty} \in \ell^{p} (\mathbb{N})$. Given two sequences $(f_{n})_{n = 1}^{\infty}$, $(g_{n})_{n = 1}^{\infty} \subseteq X$ such that $\left\| f_{n} - g_{n} \right\|_X < \frac{\varepsilon}{2^{n}}$ for each $n \in \mathbb{N}$, we have the following:
	
	\begin{displaymath}
		\left\| \sum_{n = 1}^{\infty} \alpha_{n} (f_{n} - g_{n}) \right\|_X \leq \left\| (\alpha_{n})_{n = 1}^{\infty} \right\|_{\ell^{p} (\mathbb{N})} \varepsilon.
	\end{displaymath}
	
	Consequently,
	
	\begin{displaymath}
		\left\| \sum_{n = 1}^{\infty} \alpha_{n} g_{n} \right\|_X 
		\geq \left\| \sum_{n = 1}^{\infty} \alpha_{n} f_{n} \right\|_X - \left\| (\alpha_{n})_{n = 1}^{\infty} \right\|_{\ell^{p} (\mathbb{N})} \varepsilon.
	\end{displaymath}
	
\end{lemma}

\begin{lemma} \thlabel{system of functions}  Suppose that conditions of \thref{propnon-compact} are satisfied. Let $\varepsilon > 0$. Then, there exist $k_{\varepsilon} \in \mathbb{N}$ and three sequences of functions $(f_{j})_{j = 1}^{\infty}$, $(g_{j})_{j = 1}^{\infty}$, $(h_{j})_{j = 1}^{\infty}$ such that 

\begin{enumerate}
\item For each $j \in \mathbb{N}$, $\frac{\gamma^{p}}{2^{p}} - \varepsilon < \left\| f_{j}^{\#} \right\|_{L^{q(\cdot), p} (\Omega)}^{p} $.

\item For each $j \in \mathbb{N}$, $\left\| \nabla f_{j}^{\#} \right\|_{L^{p} (\Omega)}^{p} - \varepsilon < \left\| \nabla g_{j} \right\|_{L^{p} (\Omega)} \leq \left\| \nabla f_{j}^{\#} \right\|_{L^{p} (\Omega)}^{p}=1$.

\item For each $j \in \mathbb{N}$, $\left\| f_{j}^{\#} \right\|_{L^{q(\cdot), p} (\Omega)}^{p} - \varepsilon \leq \left\| g_{j} \right\|_{L^{q(\cdot), p} (\Omega)}^{p} \leq \left\| f_{j}^{\#} \right\|_{L^{q(\cdot), p} (\Omega)}^{p}$.

\item For each $j \in \mathbb{N}$, $\left\| g_{j} \right\|_{L^{q(\cdot), p} (\Omega)}^{p} - \varepsilon \leq \left\| h_{j} \right\|_{L^{q(\cdot), p} (\Omega)}^{p} \leq \left\| g_{j} \right\|_{L^{q(\cdot), p} (\Omega)}^{p}$.

\item The functions $h_{j}$'s all have disjoint support.

\item The functions $\nabla (g_{j})$'s all have disjoint support.

\item For any $(\alpha_{j})_{j = 1}^{\infty} \in \ell^{p} (\mathbb{N})$,

\begin{equation}
\left\| \sum_{j = 1}^{\infty} \alpha_{j} h_{j} \right\|_{L^{q(\cdot), p} (\Omega)}^{p} \geq \sum_{j \in \left\{ m \in \mathbb{N} : | \alpha_{m} | > \frac{1}{2^{m + k_{\varepsilon}}} \right\} } | \alpha_{j} |^{p}  \left[ \left\| h_{j} \right\|_{L^{q(\cdot), p} (\Omega)}^{p} - \varepsilon \right].
\end{equation}

\end{enumerate}

\end{lemma}

\begin{proof}
 Without loss of generality, we can assume, in the rest of this proof, that $x_0=0$ and $\Omega = B_{a}(0)$, where $a$ is sufficiently small such that $\left| B_{a} (0) \right| \leq 1$ and $a < \eta_{0}$ from \thref{propnon-compact}.
	
\begin{displaymath}
    \text{For notational convenience, we set } \tilde{q} (x) := \frac{d p}{d - p} - \frac{C_{0}}{\left| \log \left( \frac{1}{\left| x \right|}\right)\right|}.
\end{displaymath}	
 	
Obviously, there exists $0 < m \leq a$ such that for all $r < m$, $\gamma^{p} \leq \gamma_{r}^{p} < \gamma^{p} + \varepsilon$.

Now, we will construct three sequences of functions $(f_{n})_{n = 1}^{\infty}$, $(g_{n})_{n = 1}^{\infty}$ , and $(h_{n})_{n = 1}^{\infty}$.

\bigskip

Choose $r_{1} < m$. Then, there exists $f_{1} : B_{r_{1}} (0) \to \mathbb{R}$ such that 

\begin{displaymath}
    \gamma^{p} - 2^{p} \epsilon \leq \gamma_{r_{1}}^{p} - 2^{p} \epsilon < \frac{\left\| f_{1} \right\|_{L^{\tilde{q} (\cdot) , p} \left( B_{r_{1}} (0) \right)}^{p} }{\left\| \nabla f_{1} \right\|_{L^{p} \left( B_{r_{1}} (0) \right)}^{p} } \leq \gamma_{r_{1}}^{p}.
\end{displaymath}

By \thref{symmetrization}, we have that
\begin{align*}
    & \frac{\left\| f_{1} \right\|_{L^{\tilde{q} (\cdot), p} \left( B_{r_{1}} (0) \right)}^{p} }{\left\| \nabla f_{1} \right\|_{L^{p} \left( B_{r_{1}} (0) \right)}^{p} } \leq \frac{\left\| f_{1}^{\#} \right\|_{L^{\tilde{q} (\cdot), p} \left( B_{r_{1}} (0) \right)}^{p} }{\left\| \nabla f_{1} \right\|_{L^{p} \left( B_{r_{1}} (0) \right)}^{p} } \text{ by \thref{symmetrization}} \\
    & \leq \frac{\left\| f_{1}^{\#} \right\|_{L^{\tilde{q} (\cdot), p} \left( B_{r_{1}} (0) \right)}^{p} }{\left\| \nabla f_{1}^{\#} \right\|_{L^{p} \left( B_{r_{1}} (0) \right)}^{p} } \text{ by Polya-Szego inequality} \\
    & \leq \left[  1 + \left| B_{r_{1}} (0) \right| \right]^{p} \frac{\left\| f_{1}^{\#} \right\|_{L^{q(\cdot), p} \left( B_{r_{1}} (0) \right)}^{p} }{ \left\| \nabla f_{1}^{\#} \right\|_{L^{p} \left( B_{r_{1}} (0) \right)}^{p} } \leq 2^{p} \frac{\left\| f_{1}^{\#} \right\|_{L^{q(\cdot), p} \left( B_{r_{1}} (0) \right)}^{p} }{ \left\| \nabla f_{1}^{\#} \right\|_{L^{p} \left( B_{r_{1}} (0) \right)}^{p} }.
\end{align*}

Combining the above two inequalities, we see that

\begin{displaymath}
    \frac{\gamma^{p}}{2^{p}} - \epsilon \leq  \frac{\left\| f_{1}^{\#} \right\|_{L^{q(\cdot), p} \left( B_{r_{1}} (0) \right)}^{p} }{ \left\| \nabla f_{1}^{\#} \right\|_{L^{p} \left( B_{r_{1}} (0) \right)}^{p} }.
\end{displaymath}

Thus, there exists $f_{1} : B_{r_{1} (0)} \to \mathbb{R}$ such that $\| \nabla f_{1}^{\#} \|_{L^{p} \left( B_{r_{1}} (0) \right)} = 1$ and $\frac{\gamma^{p}}{2^{p}} - \varepsilon < \| f_{1}^{\#} \|_{L^{q(\cdot), p} \left( B_{r_{1}} (0) \right)}^{p}$.

\bigskip

By the absolute continuity of both $\| \cdot \|_{L^{p}}$ and $\| \cdot \|_{L^{q( \cdot) , p}}$, there exists $w_{1} > 0$ such that $\| \nabla f_{1}^{\#} \|_{L^{p} \left( B_{w_{1}} (0) \right)}^{p} < \varepsilon$ and $\| f_{1}^{\#} \|_{L^{q(\cdot), p} \left( B_{w_{1}} (0) \right)}^{p} < \varepsilon$.

\bigskip

Since $f_{1}^{\#}$ is radially symmetric and decreasing about the point $x = 0$, the  points of discontinuity of $f_{1}^{\#}$ are at most countable. Choose $s_{1} < w_{1}$ such that $f_{1}^{\#}$ is continuous on $\{ x \in \mathbb{R}^{d} : |x| = s_{1} \}$. Choose $y_{1} \in \{ x \in \mathbb{R}^{d} : |x| = s_{1} \}$ and define $R_{1} := f_{1}^{\#} (y_{1})$.

\begin{displaymath}
\text{Define } g_{1} : B_{r_{1}} (0) \to \mathbb{R} \text{ by } g_{1} (x) :=
\begin{cases}
f_{1}^{\#} (x)
& \text{if } |x| \geq B_{s_{1}} (0),\\
R_{1} & \text{if } |x| < B_{s_{1}} (0). 
\end{cases} 
\end{displaymath}

By applying integration by parts, we obtain:

\begin{displaymath}
(\nabla g_{1}) (x) :=
\begin{cases}
(\nabla f_{1}^{\#}) (x)
& \text{if } |x| \geq B_{s_{1}} (0),\\
0 & \text{if } |x| < B_{s_{1}} (0). 
\end{cases} 
\end{displaymath}

So, we obtain the following inequalities:

$\left\| \nabla f_{1}^{\#} \right\|_{L^{p} \left( B_{r_{1}} (0) \right)}^{p} - \varepsilon < \left\| \nabla g_{1} \right\|_{L^{p} \left( B_{r_{1}} (0) \right)}^{p} \leq \left\| \nabla f_{1}^{\#} \right\|_{L^{p} \left( B_{r_{1}} (0) \right)}^{p}$ and

\bigskip

$\left\| f_{1}^{\#} \right\|_{L^{q( \cdot), p} \left( B_{r_{1}} (0) \right)}^{p} - \varepsilon < \left\| g_{1} \right\|_{L^{q(\cdot), p} \left( B_{r_{1}} (0) \right)}^{p} \leq \left\| f_{1}^{\#} \right\|_{L^{q( \cdot), p} \left( B_{r_{1}} (0) \right)}^{p}$.

\bigskip

Next, since $g_{1}$ is symmetric and radially decreasing about the point $x = 0$, by the absolute continuity of $\| \cdot \|_{L^{q(\cdot), p}}$, there exists $0 < t_{1} < s_{1}$ such that $\| g_{1} \|_{L^{q(\cdot), p} \left( B_{r_{1}} (0) \right)}^{p} < \varepsilon$.

\bigskip

Define $h_{1} : B_{r_{1}} (0) \to \mathbb{R}$ by $h_{1} := [ 1 - \chi_{B_{t_{1}} (0)} ] g_{1} $. With this, we obtain the following inequality:

$\| g_{1} \|_{L^{q(\cdot), p} \left( B_{r_{1}} (0) \right)}^{p} - \varepsilon < \| h_{1} \|_{L^{q(\cdot), p} \left( B_{r_{1}} (0) \right)}^{p} \leq \| g_{1} \|_{L^{q(\cdot), p} \left( B_{r_{1}} (0) \right)}^{p}$.

\bigskip

Next, let $\nu_{d}$ be the volume of the unit ball in $\mathbb{R}^{d}$ and choose $\delta_{2} > 0$ such that $\frac{1}{p} (\nu_{d} \delta_{2}^{n} )^{\frac{p}{q_{+}}} (1 + \nu_{d} \delta_{2}^{n} ) f_{1}(y) 2^{2 + k_{\varepsilon}} < \varepsilon$. Define $r_{2} := \min \{ t_{1}, \delta_{2} \}$. There exists $f_{2} : B_{r_{2}} (0) \to \mathbb{R}$ such that 

\begin{displaymath}
    \gamma^{p} - 2^{p} \epsilon \leq \gamma_{r_{2}}^{p} - 2^{p} \epsilon < \frac{\left\| f_{2} \right\|_{L^{\tilde{q} (\cdot), p} \left( B_{r_{2}} (0) \right)}^{p} }{\left\| \nabla f_{2} \right\|_{L^{p} \left( B_{r_{2}} (0) \right)}^{p} } \leq \gamma_{r_{2}}^{p}.
\end{displaymath}

Then, we repeat the same process inductively as above to obtain the four sequences of functions $(f_{n})_{n = 1}^{\infty}$, $\left( f_{n}^{\#} \right)_{n = 1}^{\infty}$, $(g_{n})_{n= 1}^{\infty}$ , and $(h_{n})_{n = 1}^{\infty}$, where $f_{n}$, $f_{n}^{\#}$, $g_{n}$, and $h_{n}$ are functions supported on $B_{r_{n}} (0)$.

\bigskip

By the construction of both $g_{n}$ and $h_{n}$, we see that $\sup_{x \in B_{r_{n}} (0)} | g_{n} (x) |= R_{n} = \sup_{x \in B_{r_{n}} (0)} | h_{n} (x) |$. 

\bigskip

Let $(\alpha_{n})_{n = 1}^{\infty} \in \ell^{p}(\mathbb{N})$. By \thref{main}, there exists $k \in \mathbb{N}$ ($k$ depends on $\varepsilon$ and $\left\| (\alpha_{n})_{n = 1}^{\infty} \right\|_{\ell^{p}})$ such that   

\begin{displaymath}
    \sum_{j \in \left\{ m \in \mathbb{N} : | \alpha_{m} | \leq \frac{1}{2^{m + k}} \right\}} | \alpha_{j} |^{p} < \varepsilon \sum_{j = 1}^{\infty} | \alpha_{j} |^{p} 
\end{displaymath}

Observe that for each $n \in \left\{ m \in \mathbb{N} : | \alpha_{m} | > \frac{1}{2^{m + k_{\varepsilon}}} \right\}$, we have

\begin{align*}
    & \int_{0}^{\frac{R_{n - 1}}{|\alpha_{n}|}} \lambda^{p - 1} \left\| \chi_{\left\{ x \in B_{r_{n}} (0) : | h_{n} (x) | > \lambda \right\}} \right\|_{L^{q(\cdot)} [B_{r_{n}} (0)]}^{p} d \lambda \\
    & \leq  \int_{0}^{\frac{R_{n - 1}}{|\alpha_{n}|}} \lambda^{p - 1} \left\| \chi_{\left\{ x \in B_{r_{n}} (0) : | h_{n} (x) | > \lambda \right\}} \right\|_{L^{q_{+}} [B_{r_{n}} (0)]}^{p} (1 + \nu_{d} r_{n}^{d})d \lambda \\
    & \leq  \int_{0}^{\frac{R_{n - 1}}{|\alpha_{n}|}} \lambda^{p - 1} (\nu_{d} r_{n}^{d})^{\frac{p}{q_{+}}} (1 + \nu_{d} r_{n}^{d}) d \lambda \\
    & = \frac{1}{p} (\nu_{d} r_{n}^{d})^{\frac{p}{q_{+}}} (1 + \nu_{d} r_{n}^{d}) \left( \frac{R_{n - 1}}{|\alpha_{n}|} \right)^{p} \\
    & \leq \frac{1}{p} (\nu_{d} r_{n}^{d})^{\frac{p}{q_{+}}} 
    (1 + \nu_{d} r_{n}^{d}) \left( R_{n - 1} 2^{n + k_{\varepsilon}} \right)^{p} < \varepsilon.
\end{align*}

Finally, we proceed to prove (7). 





\begin{align*}
    & \left\| \sum_{j = 1}^{\infty} \alpha_{j} h_{j} \right\|_{L^{q(\cdot), p} (\Omega)}^{p} = \int_{0}^{\infty} \lambda^{p - 1} \left\| \chi_{\left\{ x \in \Omega : \left| \sum_{j = 1}^{\infty} \alpha_{j} h_{j} (x) \right| > \lambda \right\}} \right\|_{L^{q(\cdot)} (\Omega)}^{p} d \lambda \\
    & = \int_{0}^{\infty} \lambda^{p - 1} \left\| \sum_{j = 1}^{\infty} \chi_{\left\{ x \in \Omega : \left| \alpha_{j} h_{j} (x) \right| > \lambda \right\}} \right\|_{L^{q(\cdot)} (\Omega)}^{p} d \lambda \quad \text{ (since each } h_{j} \text{ has disjoint support)} \\
    & \geq \int_{0}^{\infty} \lambda^{p - 1} \left\| \sum_{j \in \left\{ m \in \mathbb{N} : | \alpha_{m} | > \frac{1}{2^{m + k_{\varepsilon}}} \right\}} \chi_{\left\{ x \in \Omega : \left| \alpha_{j} h_{j} (x) \right| > \lambda \right\}} \right\|_{L^{q(\cdot)} (\Omega)}^{p} d \lambda \\
    & \geq \int_{0}^{\infty} \lambda^{p - 1} \left\| \sum_{j \in \left\{ m \in \mathbb{N} : | \alpha_{m} | > \frac{1}{2^{m + k_{\varepsilon}}} \right\}} \chi_{\left\{ x \in \Omega : \left| \alpha_{j} h_{j} (x) \right| > \lambda \right\}} \chi_{[R_{j - 1}, R_{j} )} (\lambda) \right\|_{L^{q(\cdot)} (\Omega)}^{p} d \lambda \\
    & = \int_{0}^{\infty} \lambda^{p - 1} \sum_{j \in \left\{ m \in \mathbb{N} : | \alpha_{m} | > \frac{1}{2^{m + k_{\varepsilon}}} \right\}} \left\| \chi_{\left\{ x \in \Omega : \left| \alpha_{j} h_{j} (x) \right| > \lambda \right\}} \chi_{[R_{j - 1}, R_{j} )} (\lambda) \right\|_{L^{q(\cdot)} (\Omega)}^{p} d \lambda \\
    & = \sum_{j \in \left\{ m \in \mathbb{N} : | \alpha_{m} | > \frac{1}{2^{m + k_{\varepsilon}}} \right\}} \int_{R_{j - 1}}^{R_{j}} \lambda^{p - 1} \left\| \chi_{\left\{ x \in \Omega : \left| \alpha_{j} h_{j} (x) \right| > \lambda \right\}} \right\|_{L^{q(\cdot)} (\Omega)}^{p} d \lambda \\
    & = \sum_{j \in \left\{ m \in \mathbb{N} : | \alpha_{m} | > \frac{1}{2^{m + k_{\varepsilon}}} \right\}} |\alpha_{j}|^{p} \int_{\frac{R_{j - 1}}{|\alpha_{j}|}}^{\frac{R_{j}}{|\alpha_{j}|}} \lambda^{p - 1} \left\| \chi_{\left\{ x \in \Omega : \left| h_{j} (x) \right| > \lambda \right\}} \right\|_{L^{q(\cdot)} (\Omega)}^{p} d \lambda \\
    & = \sum_{j \in \left\{ m \in \mathbb{N} : | \alpha_{m} | > \frac{1}{2^{m + k_{\varepsilon}}} \right\}} |\alpha_{j}|^{p} \int_{\frac{R_{j - 1}}{|\alpha_{j}|}}^{R_{j}} \lambda^{p - 1} \left\| \chi_{\left\{ x \in \Omega : \left| h_{j} (x) \right| > \lambda \right\}} \right\|_{L^{q(\cdot)}(\Omega)}^{p} d \lambda \\
    & \geq \sum_{j \in \left\{ m \in \mathbb{N} : | \alpha_{m} | > \frac{1}{2^{m + k_{\varepsilon}}} \right\}} |\alpha_{j}|^{p} \left[ \int_{0}^{R_{j}} \lambda^{p - 1} \left\| \chi_{\left\{ x \in \Omega : \left| h_{j} (x) \right| > \lambda \right\}} \right\|_{L^{q(\cdot)} (\Omega)}^{p} d \lambda - \varepsilon \right] \\
    & = \sum_{j \in \left\{ m \in \mathbb{N} : | \alpha_{m} | > \frac{1}{2^{m + k_{\varepsilon}}} \right\}} |\alpha_{j}|^{p} \left[ \left\| h_{j} \right\|_{L^{q(\cdot), p} (\Omega)}^{p} - \varepsilon \right].
\end{align*}

\end{proof}

With the above construction of the sequences $(f_{n})_{n = 1}^{\infty}$, $\left( f_{n}^{\#} \right)_{n = 1}^{\infty}$, $(g_{n})_{n = 1}^{\infty}$ and $(h_{n})_{n = 1}^{\infty}$, we can show that there is  a lower bound for the Bernstein numbers. 

\begin{theorem} \thlabel{mainresult}
    Consider the embedding $I : W_{0}^{1, p} (\Omega) \to L^{q(\cdot), p} (\Omega)$ with conditions from  \thref{propnon-compact}. Then, for each $N \in \mathbb{N}$, $b_{N} (I) \geq \frac{\gamma}{2} > 0$, where $\gamma$ is defined as in \eqref{gamma}.
\end{theorem}

\begin{proof}
    
Let $N \in \mathbb{N}$. Let $V_{\epsilon, N}$ be the $n$-dimensional subspace of $W_0^{1, p}$ spanned by $g_{\epsilon, 1}, \cdots, g_{\epsilon, N}$ as constructed in \thref{system of functions}. By the definition of the Bernstein number, for all $\epsilon > 0$,

\begin{displaymath}
    b_{N} (I) \geq \inf_{x \in V_{\epsilon, N}, \| x \|_{W_0^{1, p} } = 1} \left\| x \right\|_{L^{q(\cdot), p}}
\end{displaymath}

Now, fix $\epsilon > 0$ and consider $V_{\epsilon, N}$. Let $x = \sum_{j = 1}^{N} \alpha_{j} g_{j}$ be an arbitrary element of $V_{\epsilon, N}$, with $\sum_{j=1}^N |\alpha_j |^p=1$. Applying 
\thref{almost orthogonal}, \thref{system of functions} and \thref{main}, we have:
\begin{align*}
    & \frac{\left\| \sum_{j = 1}^{N} \alpha_{j} g_{j} \right\|_{L^{q(\cdot), p}}}{\left\| \sum_{j = 1}^{N} \alpha_{j} \nabla g_{j} \right\|_{L^{p}}} \geq \frac{\left\| \sum_{j = 1}^{N} \alpha_{j} h_{j} \right\|_{L^{q(\cdot), p}} - \left\| (\alpha_{j} )_{j = 1}^{N} \right\|_{\ell^{\infty} (\mathbb{N})}  \varepsilon }{\left( \sum_{j = 1}^{N} | \alpha_{j} |^{p} \left\| \nabla g_{j} \right\|_{L^{p}}^{p} \right)^{\frac{1}{p}}} \text{ by \thref{almost orthogonal}}\\
    & \geq \frac{\left[ \sum_{j \in \left\{ m \in \mathbb{N} : | \alpha_{m} | > \frac{1}{2^{m + k_{\varepsilon}}} \right\}} |\alpha_{j}|^{p} \left( \left\| h_{j} \right\|_{L^{q(\cdot), p}}^{p} - \varepsilon \right) \right]^{\frac{1}{p}} - \left\| (\alpha_{j} )_{j = 1}^{N} \right\|_{\ell^{\infty} (\mathbb{N})}  \varepsilon }{\left( \sum_{j = 1}^{N} | \alpha_{j} |^{p} \left\| \nabla g_{j} \right\|_{L^{p}}^{p} \right)^{\frac{1}{p}}} \text{by \thref{system of functions} (7)}\\
    & \geq \frac{\left[ \sum_{j \in \left\{ m \in \mathbb{N} : | \alpha_{m} | > \frac{1}{2^{m + k_{\varepsilon}}} \right\}} |\alpha_{j}|^{p} \left( \left\| g_{j} \right\|_{L^{q(\cdot), p}}^{p} - 2 \varepsilon \right) \right]^{\frac{1}{p}} - \left\| (\alpha_{j} )_{j = 1}^{N} \right\|_{\ell^{\infty} (\mathbb{N})}  \varepsilon }{\left( \sum_{j = 1}^{N} | \alpha_{j} |^{p} \left\| \nabla g_{j} \right\|_{L^{p}}^{p} \right)^{\frac{1}{p}}} \text{by \thref{system of functions} (4)} \\
    & \geq \frac{\left[ \sum_{j \in \left\{ m \in \mathbb{N} : | \alpha_{m} | > \frac{1}{2^{m + k_{\varepsilon}}} \right\}} |\alpha_{j}|^{p} \left( \left\| f_{j}^{\#} \right\|_{L^{q(\cdot), p}}^{p} - 3 \varepsilon \right) \right]^{\frac{1}{p}} - \left\| (\alpha_{j} )_{j = 1}^{N} \right\|_{\ell^{\infty} (\mathbb{N})}  \varepsilon }{\left( \sum_{j = 1}^{N} | \alpha_{j} |^{p} \left\| \nabla f_{j} \right\|_{L^{p}}^{p} \right)^{\frac{1}{p}}} \text{by \thref{system of functions} (3)} \\
    & \geq \frac{\left[ \sum_{j \in \left\{ m \in \mathbb{N} : | \alpha_{m} | > \frac{1}{2^{m + k_{\varepsilon}}} \right\}} |\alpha_{j}|^{p} \left( \frac{\gamma^{p}}{2^{p}} - 4 \varepsilon \right) \right]^{\frac{1}{p}} - \left\| (\alpha_{j} )_{j = 1}^{N} \right\|_{\ell^{\infty} (\mathbb{N})}  \varepsilon }{\left( \sum_{j = 1}^{N} | \alpha_{j} |^{p} \right)^{\frac{1}{p}}} \text{by \thref{system of functions} (1)} \\
    & = \frac{\left( \frac{\gamma^{p}}{2^{p}} - 4 \varepsilon \right)^{\frac{1}{p}} \left[ \sum_{j \in \left\{ m \in \mathbb{N} : | \alpha_{m} | > \frac{1}{2^{m + k_{\varepsilon}}} \right\}} |\alpha_{j}|^{p} \right]^{\frac{1}{p}} - \left\| (\alpha_{j} )_{j = 1}^{N} \right\|_{\ell^{\infty} (\mathbb{N})}  \varepsilon }{\left( \sum_{j = 1}^{N} | \alpha_{j} |^{p} \right)^{\frac{1}{p}}}  \\
    & \geq \frac{\left( \frac{\gamma^{p}}{2^{p}} - 4 \varepsilon \right)^{\frac{1}{p}} \left[ (1 - \varepsilon) \sum_{j = 1}^{N} |\alpha_{j}|^{p} \right]^{\frac{1}{p}} - \left\| (\alpha_{j} )_{j = 1}^{N} \right\|_{\ell^{\infty} (\mathbb{N})}  \varepsilon }{\left( \sum_{j = 1}^{N} | \alpha_{j} |^{p} \right)^{\frac{1}{p}}} \text{ by \thref{main} }\\
    & = \left( \frac{\gamma^{p}}{2^{p}} - 4 \varepsilon \right)^{\frac{1}{p}} (1 - \varepsilon)^{\frac{1}{p}} - \frac{\left\| (\alpha_{j} )_{j = 1}^{N} \right\|_{\ell^{\infty} (\mathbb{N})}  \varepsilon }{\left( \sum_{j = 1}^{N} | \alpha_{j} |^{p} \right)^{\frac{1}{p}}} \\
    & \geq \left( \frac{\gamma^{p}}{2^{p}} - 4 \varepsilon \right)^{\frac{1}{p}} (1 - \varepsilon)^{\frac{1}{p}} - \varepsilon.
\end{align*}

Since $x \in V_{\epsilon, N}$ is arbitrary, we obtain that:

\begin{displaymath}
    b_{N} (I) \geq \inf_{x \in V_{\epsilon, N}, \| x \|_{W^{1, p} (\Omega)} = 1} \left\| x \right\|_{L^{q(\cdot), p}} \geq \left( \frac{\gamma^{p}}{2^{p}} - 4 \varepsilon \right)^{\frac{1}{p}} (1 - \varepsilon)^{\frac{1}{p}} - \varepsilon.
\end{displaymath}

Since $\varepsilon > 0$ is arbitrary, $b_{N} (I)  \geq \frac{\gamma}{2} > 0,$ where the last inequality follows from \thref{gamma lemma}.

\end{proof}

\begin{remark}
    It follows that $\beta (I) \geq \frac{\gamma}{2} > 0$.
\end{remark}

\begin{theorem} \thlabel{mainTh} 
    Let $1 \leq p < d$, where $d \geq 2$ is a positive integer. Let $\Omega$ be a bounded domain on $\mathbb{R}^{d}$. Let $q : \Omega \to [1, \infty)$ be a measurable function satisfying $1 \leq q(x) \leq \frac{d p}{d - p}$ for almost every $x \in \Omega$. Suppose that there exist $x_{0} \in \Omega$ and constants $C_{0}$, $\eta_{0} > 0$ such that 
	
	\begin{equation} \label{non-compact3}
		q(x) \geq \frac{d p}{d - p} - \frac{C_{0}}{ \left| \log \left( \frac{1}{| x - x_{0} |} \right) \right|} \text{ for almost every } x \in \Omega \text{ with } | x - x_{0} | < \eta_{0}.
	\end{equation}
	
	Then the embedding $I : W_{0}^{1, p} (\Omega) \to L^{q(\cdot), p} (\Omega)$ is not strictly singular.
\end{theorem}

\begin{proof}
    We need to show that there exists an infinite dimensional, closed subspace $Z$ of $W_{0}^{1, p} (\Omega)$ such that 

    \begin{displaymath}
        \inf \left\{ \left\| I (x) \right\|_{L^{q(\cdot), p} (\Omega)} : \left\| \nabla x \right\|_{L^{p} (\Omega)} = 1, x \in Z \right\} > 0.
    \end{displaymath}

    First, choose $\varepsilon > 0$ such that $\left( \frac{\gamma^{p}}{2^{p}} - 4 \varepsilon \right)^{\frac{1}{p}} (1 - \varepsilon)^{\frac{1}{p}} - \varepsilon > 0$. Then, consider the infinite dimensional subspace $Z_{\epsilon} \subseteq W_{0}^{1, p} (\Omega)$ defined by: 

    \begin{displaymath}
        Z_{\varepsilon} := \overline{\text{Span} \left\{ g_{\epsilon, j} \right\}_{j = 1}^{\infty} }^{W_{0}^{1, p} (\Omega)}, \text{ where } g_{\epsilon, j} \text{ is constructed in \thref{system of functions}}.
    \end{displaymath}

    By the proof of \thref{mainresult}, for each $x \in Z_{\epsilon}$ of the form $x = \sum_{j = 1}^{N} \alpha_{j} g_{j}$, we have that:
    \begin{displaymath}
        \frac{\left\| \sum_{j = 1}^{N} \alpha_{j} g_{j} \right\|_{L^{q(\cdot), p} (\Omega)}}{ \left\| \sum_{j = 1}^{N} \alpha_{j} \nabla g_{j} \right\|_{L^{p} (\Omega)}} \geq \left( \frac{\gamma^{p}}{2^{p}} - 4 \varepsilon \right)^{\frac{1}{p}} (1 - \varepsilon)^{\frac{1}{p}} - \varepsilon.
    \end{displaymath}

    Now, let $y \in Z_{\varepsilon}$ where $\left\| \nabla y \right\|_{L^{p} (\Omega)} = 1$. Choose $ \delta > 0$ sufficiently small such that $\left\| I \right\|_{W_{0}^{1, p} \to L^{q(\cdot), p}} \delta < \left( \frac{\gamma^{p}}{2^{p}} - 4 \varepsilon \right)^{\frac{1}{p}} (1 - \varepsilon)^{\frac{1}{p}} - \varepsilon$. Then, there exists $x = \sum_{j = 1}^{N} \alpha_{j} g_{j}$ such that $\| \nabla x \|_{L^{p} (\Omega)} = 1$ and $\left\| \nabla y - \nabla x \right\|_{L^{p} (\Omega)} < \delta$.

    \begin{align*}
        & \left\| y \right\|_{L^{q(\cdot), p} (\Omega)} \geq \left\| x \right\|_{L^{q(\cdot), p} (\Omega)} - \left\| y - x \right\|_{L^{q(\cdot), p} (\Omega)} \geq \left\| x \right\|_{L^{q(\cdot), p} (\Omega)} - \left\| I \right\| \left\| \nabla y - \nabla x \right\|_{L^{p} (\Omega)} \\
        & \geq \left( \frac{\gamma^{p}}{2^{p}} - 4 \varepsilon \right)^{\frac{1}{p}} (1 - \varepsilon)^{\frac{1}{p}} - \varepsilon - \left\| I \right\| \delta > 0.
    \end{align*}

    Since $y \in Z_{\varepsilon}$ is arbitrary, 
    \begin{displaymath}
        \inf \left\{ \left\| I (y) \right\|_{L^{q(\cdot), p} (\Omega)} : \left\| \nabla y \right\|_{L^{p} (\Omega)} = 1, y \in Z_{\epsilon} \right\} \geq \left( \frac{\gamma^{p}}{2^{p}} - 4 \varepsilon \right)^{\frac{1}{p}} (1 - \varepsilon)^{\frac{1}{p}} - \varepsilon - \left\| I \right\| \delta > 0.
    \end{displaymath}
\end{proof}

{\bf Acknowledgements.} The authors thank the anonymous referee for reading the manuscript carefully and providing valuable comments.
\bibliography{bibliography}

\begin{thebibliography}{11}
\providecommand{\natexlab}[1]{#1}
\providecommand{\url}[1]{\texttt{#1}}
\expandafter\ifx\csname urlstyle\endcsname\relax
  \providecommand{\doi}[1]{doi: #1}\else
  \providecommand{\doi}{doi: \begingroup \urlstyle{rm}\Url}\fi

\bibitem[Bouchala(2020)]{Bo2020}
Ond\v{r}ej Bouchala.
\newblock Measures of non-compactness and {S}obolev-{L}orentz spaces.
\newblock \emph{Z. Anal. Anwend.}, 39\penalty0 (1):\penalty0 27--40, 2020.
\newblock ISSN 0232-2064.
\newblock \doi{10.4171/zaa/1649}.

\bibitem[Edmunds and Evans(1987)]{EdEv}
D.~E. Edmunds and W.~D. Evans.
\newblock \emph{Spectral theory and differential operators}.
\newblock Oxford Mathematical Monographs. The Clarendon Press, Oxford University Press, New York, 1987.
\newblock ISBN 0-19-853542-2.

\bibitem[Edmunds et~al.(2023)Edmunds, Gogatishvili, and Nekvinda]{EGN}
D.~E. Edmunds, A.~Gogatishvili, and A.~Nekvinda.
\newblock Almost compact and compact embeddings of variable exponent spaces.
\newblock \emph{Studia Math.}, 268\penalty0 (2):\penalty0 187--211, 2023.
\newblock ISSN 0039-3223.
\newblock \doi{10.4064/sm211206-24-2}.

\bibitem[Hencl(2003)]{He}
Stanislav Hencl.
\newblock Measures of non-compactness of classical embeddings of {S}obolev spaces.
\newblock \emph{Math. Nachr.}, 258:\penalty0 28--43, 2003.
\newblock ISSN 0025-584X.
\newblock \doi{10.1002/mana.200310085}.

\bibitem[Kempka and Vybíral(2014-06)]{KeVy}
Henning Kempka and Jan Vybíral.
\newblock Lorentz spaces with variable exponents: Lorentz spaces with variable exponents.
\newblock \emph{Math. Nachr.}, 287\penalty0 (8):\penalty0 938--954, 2014-06.
\newblock ISSN 0025584X.
\newblock \doi{10.1002/mana.201200278}.
\newblock URL \url{https://onlinelibrary.wiley.com/doi/10.1002/mana.201200278}.

\bibitem[Kerman and Pick(2006)]{KP2006}
Ron Kerman and Lubo\v{s} Pick.
\newblock Optimal {S}obolev imbeddings.
\newblock \emph{Forum Math.}, 18\penalty0 (4):\penalty0 535--570, 2006.
\newblock ISSN 0933-7741.
\newblock \doi{10.1515/FORUM.2006.028}.

\bibitem[Kov\'{a}\v{c}ik and R\'{a}kosn\'{\i}k(1991)]{KovRak}
Ondrej Kov\'{a}\v{c}ik and Ji\v{r}\'{\i} R\'{a}kosn\'{\i}k.
\newblock On spaces {$L^{p(x)}$} and {$W^{k,p(x)}$}.
\newblock \emph{Czechoslovak Math. J.}, 41(116)\penalty0 (4):\penalty0 592--618, 1991.
\newblock ISSN 0011-4642.

\bibitem[Kurata and Shioji(2008-03)]{KS2008}
Kazuhiro Kurata and Naoki Shioji.
\newblock Compact embedding from {$W_0^{1,2}(\Omega) \to L^{q(x)}(\Omega)$} and its application to nonlinear elliptic boundary value problem with variable critical exponent.
\newblock \emph{Journal of Mathematical Analysis and Applications}, 339\penalty0 (2):\penalty0 1386--1394, 2008-03.
\newblock ISSN 0022247X.
\newblock \doi{10.1016/j.jmaa.2007.07.083}.
\newblock URL \url{https://linkinghub.elsevier.com/retrieve/pii/S0022247X07010116}.

\bibitem[Lang and Mihula(2023)]{LaMi}
Jan Lang and Zden\v{e}k Mihula.
\newblock Different degrees of non-compactness for optimal {S}obolev embeddings.
\newblock \emph{J. Funct. Anal.}, 284\penalty0 (10):\penalty0 Paper No. 109880, 22, 2023.
\newblock ISSN 0022-1236.
\newblock \doi{10.1016/j.jfa.2023.109880}.

\bibitem[Lang et~al.(2021)Lang, Musil, Ol\v{s}\'{a}k, and Pick]{LMOP}
Jan Lang, V\'{\i}t Musil, Miroslav Ol\v{s}\'{a}k, and Lubo\v{s} Pick.
\newblock Maximal non-compactness of {S}obolev embeddings.
\newblock \emph{J. Geom. Anal.}, 31\penalty0 (9):\penalty0 9406--9431, 2021.
\newblock ISSN 1050-6926.
\newblock \doi{10.1007/s12220-020-00522-y}.

\bibitem[Peetre(1966)]{Pe}
Jaak Peetre.
\newblock Espaces d'interpolation et th\'{e}or\`eme de {S}oboleff.
\newblock \emph{Ann. Inst. Fourier (Grenoble)}, 16\penalty0 (fasc. 1):\penalty0 279--317, 1966.
\newblock ISSN 0373-0956.

\end{thebibliography}
\end{document}